\documentclass[aos,MSNbibl,seceqn,dvips]{arximspdf}

\doi{10.1214/13-AOS1088} 
\volume{41}
\issue{2}
\pubyear{2013}
\firstpage{604}
\lastpage{640}

\makeatletter

\newcommand{\rrVert}{\Vert}
\newcommand{\llVert}{\Vert}
\def\cal{\mathcal}
\def\eps{\varepsilon}

\newtheorem{proposition}{Proposition}
\newtheorem{lemma}{Lemma}
\newtheorem{theorem}{Theorem}
\newtheorem{corollary}{Corollary}
\newproclaim{remark}{Remark}
\newproclaim{remarks}{Remarks}
\newproclaim{example}{Example}

\def\EE{\mathbb{E}}

\makeatother

\begin{document}
\begin{frontmatter}

\title{Low rank estimation of smooth kernels on graphs}
\runtitle{Low rank estimation of smooth kernels}

\begin{aug}
\author{\fnms{Vladimir} \snm{Koltchinskii}\corref{}\ead[label=e1]{vlad@math.gatech.edu}\thanksref{t1}}
\and
\author{\fnms{Pedro} \snm{Rangel}\ead[label=e2]{prangel@math.gatech.edu}\thanksref{t2}}
\thankstext{t1}{Supported in part by NSF Grants DMS-12-07808,
DMS-09-06880 and CCF-0808863.}
\thankstext{t2}{Supported in part by NSF Grant CCF-0808863.}
\runauthor{V. Koltchinskii and P. Rangel}
\affiliation{Georgia Institute of Technology}
\address{School of Mathematics\\
Georgia Institute of Technology\\
Atlanta, Georgia 30332-0160\\
USA\\
\printead{e1}\\
\phantom{E-mail:\ }\printead*{e2}}
\end{aug}

\received{\smonth{7} \syear{2012}}
\revised{\smonth{1} \syear{2013}}

%
\begin{abstract}
Let $(V,A)$ be a weighted graph with a finite vertex set $V$, with a
symmetric matrix of
nonnegative weights $A$ and with Laplacian $\Delta$. Let $S_{\ast
}\dvtx V\times V\mapsto{\mathbb R}$ be a symmetric kernel defined on the
vertex set $V$.
Consider $n$ i.i.d. observations $(X_j,X_j', Y_j), j=1,\ldots, n$, where
$X_j, X_j'$ are independent random vertices sampled from the uniform
distribution
in $V$ and $Y_j\in{\mathbb R}$ is a real valued response variable such that
${\mathbb E}(Y_j|X_j,X_j')=S_{\ast}(X_j,X_j'), j=1,\ldots, n$.
The goal is to estimate the kernel $S_{\ast}$ based on the data
$(X_1,X_1', Y_1),\ldots, (X_n, X_n', Y_n)$ and under the assumption
that $S_{\ast}$ is low rank and, at the same time, smooth
on the graph (the smoothness being characterized by discrete Sobolev norms
defined in terms of the graph Laplacian). We obtain several
results for such problems including minimax lower bounds on the $L_2$-error
and upper bounds for penalized least squares estimators both with
nonconvex and with convex penalties.
\end{abstract}

%


\end{frontmatter}

\section{Introduction}\label{intro}

We study a problem of estimation of a symmetric
kernel $S_{\ast}\dvtx V\times V\mapsto{\mathbb R}$ defined on a large weighted
graph with a vertex set $V$ and $m:=\operatorname{ card}(V)$ based on a finite
number of noisy linear measurements of $S_{\ast}$.
For simplicity, assume that these are the measurements of randomly picked
entries of $m\times m$ matrix $(S_{\ast}(u,v))_{u,v\in V}$, which is a
standard sampling model in matrix completion. More precisely, let
$(X_j, X_j',Y_j), j=1,\ldots, n$ be $n$ independent copies of a random
triple $(X,X',Y)$, where $X,X'$ are independent random vertices sampled
from the uniform distribution $\Pi$ in $V$, and $Y\in{\mathbb R}$ is
a ``measurement'' of the kernel $S_{\ast}$ at a random location
$(X,X')$ in the sense that ${\mathbb E}(Y|X,X')=S_{\ast}(X,X')$.
In what follows, we assume that, for some constant $a>0$, $|Y|\leq a$ a.s.,
which implies that $|S_{\ast}(u,v)|\leq a, u,v\in V$.
The target kernel $S_{\ast}$ is to be estimated based on its i.i.d.
measurements $(X_j, X_j',Y_j), j=1,\ldots, n$. We would like to study
this problem in the case when the target kernel $S_{\ast}$ is, on the
one hand, ``low rank'' [i.e., $\operatorname{ rank}(S_{\ast})$ is relatively
small compared to $m$], and on the other hand, it is ``smooth'' in the
sense that its ``Sobolev-type norm'' is not too large.
Discrete versions of Sobolev norms can be defined for functions and kernels
on weighted graphs in terms of their graph Laplacians.
The problem of estimation of smooth low-rank kernels is of importance
in a number of applications,
such as learning kernels representing and predicting similarities between
objects, various classification problems in large complex networks
(e.g., edge sign prediction) as well as matrix completion problems
in the design of recommender systems (collaborative filtering).
Our main motivation, however, is mostly
theoretical: we would like to explore
to which extent taking into account smoothness of the target kernel
could improve the existing methods of low rank recovery.

We introduce some notation used throughout the paper.
Let $\mathcal{S}_V$ be the linear space of \emph{symmetric kernels}
$S\dvtx V\times V\mapsto{\mathbb R}$, $S(u,v)=S(v,u), u,v\in V$ (or,
equivalently, symmetric $m\times m$ matrices with real entries). Given
$S \in\mathcal{S}_V$,
we use the notation $\operatorname{ rank}(S)$ for the rank of $S$ and $\operatorname{ tr}(S)$
for its trace. For two functions $f,g\dvtx V\mapsto{\mathbb R}$, $(f\otimes
g)(u,v):=f(u)g(v)$. Suppose that $S=\sum_{j=1}^{r}{\mu_j (\psi_j
\otimes\psi_j)}$
is the spectral representation of $S$ with $r=\operatorname{ rank}(S)$,
$\mu_1, \ldots, \mu_r$ being nonzero eigenvalues of $S$ repeated
with their multiplicities and $\psi_1,\ldots,\psi_r$ being the
corresponding orthonormal eigenfunctions (obviously, there are multiple
choices of $\psi_j$s in the case of repeated eigenvalues).
We will define $\operatorname{ sign}(S)$ as
$\operatorname{ sign}(S):=\sum_{j=1}^{r}{\operatorname{ sign}(\mu_j)(\psi_j \otimes\psi
_j})$ and
the support of $S$ as $\operatorname{ supp}(S):=\mbox{l.s.}\{ \psi_1,\ldots,\psi
_r\}$.\setcounter{footnote}{2}\footnote{``l.s.'' means ``the linear span.''}
For $1\leq p< \infty$, define the Schatten $p$-norm of $S$ as
$
\|S\|_p:=(\operatorname{ tr}(|S|^p))^{1/p}= ( \sum_{j=1}^{r} |\mu_j|^p
)^{1/p},
$
where $|S|:=\sqrt{S^2}$.
For $p=1$, $\|\cdot\|_1$ is also called the nuclear norm and, for
$p=2$, $\|\cdot\|_2$ is called the Hilbert--Schmidt or Frobenius norm.
This norm is induced by the Hilbert--Schmidt inner product which will
be denoted by $\langle\cdot,\cdot\rangle$. The operator norm of $S$
is defined as $\|S\|:=\max_{j}|\mu_j|$.\footnote{With some abuse of
notation, we also denote occasionally the canonical Euclidean inner
product in ${\mathbb R}^V$ by $\langle\cdot,\cdot\rangle$ and
the corresponding Euclidean norm by $\|\cdot\|$.}

Let $\Pi^2:=\Pi\otimes\Pi$ be the distribution of random
couple $(X,X')$. The $L_2(\Pi^2)$-norm of kernel $S$,
\[
\|S\|_{L_2(\Pi^2)}^2=\int_{V\times V}\bigl|S(u,v)\bigr|^2
\Pi^2(du,dv)={\mathbb E}\bigl|S\bigl(X,X'
\bigr)\bigr|^2,
\]
is naturally related to the sampling model studied in the paper, and it
will be used
to measure the estimation error.
Denote by
$\langle\cdot,\cdot\rangle_{L_2(\Pi^2)}$ the corresponding inner product.
Since $\Pi$ is the uniform distribution in $V$,
$\|S\|_{L_2(\Pi^2)}^2=m^{-2}\|S\|_2^2$ and
$\langle S_1,S_2 \rangle_{L_2(\Pi^2)}=m^{-2}\langle S_1,S_2\rangle$.
In what follows, it will be often more convenient to use these rescaled versions
rather than the actual Hilbert--Schmidt norm or inner product.

We will denote by $\{e_v\dvtx v\in V\}$ the canonical orthonormal basis
of the space ${\mathbb R}^V$. Based on this basis, one can construct
matrices $E_{u,v}=E_{v,u}=\frac{1}{2}(e_u\otimes e_v+e_v\otimes e_u)$.
If $v_1,\ldots, v_m$ is an arbitrary ordering of the vertices in $V$,
then $\{E_{v_j,v_j}\dvtx j=1,\ldots, m\}\cup\{\sqrt{2}E_{v_i,v_j}\dvtx 1\leq
i<j\leq m\}$
is an orthonormal basis of the space ${\cal S}_V$ of symmetric matrices
with Hilbert--Schmidt inner product.

In standard matrix completion problems, $V$ is a finite set with no further
structure (i.e., the set of edges of the graph or the weight matrix
are not specified). In the noiseless matrix completion problems, the
target matrix $S_{\ast}$ is to be recovered from the measurements
$(X_j, X_j',Y_j), j=1,\ldots, n$,
where $Y_j=S_{\ast}(X_j,X_j')$. The following method is based on
nuclear norm minimization over the space of all matrices that ``agree''
with the data
%
\begin{equation}
\label{nucemin} \hat S:=\operatorname{ argmin}\bigl\{\|S\|_1\dvtx S
\in{\cal S}_V, S\bigl(X_j,X_j'
\bigr)=Y_j, j=1,\ldots, n\bigr\},
\end{equation}
It has been studied in detail in the recent literature; see
\cite{CandesRecht,Recht1,CandesTao,Gross-2} and references
therein. Clearly, there are low rank matrices $S_{\ast}$
that cannot be recovered based on a random sample of $n$ entries unless
$n$ is comparable with the total number of the entries of the matrix.
For instance, for given $u,v\in V$, let $S_{\ast}=E_{u,v}$.
Then, $\operatorname{ rank}(S_{\ast})\leq2$. However,
the probability that the only two nonzero entries of $S_{\ast}$
are not present in the sample is $(1-\frac{2}{m^2})^n$, and it is
close to $1$ when $n=o(m^2)$. In this case, the matrix $S_{\ast}$
cannot be recovered. So-called \emph{low coherence} assumptions have
been developed
to define classes of ``generic'' matrices that
are not ``low rank'' and ``sparse'' at the same time and for which noiseless
low rank recovery is possible with a relatively small number of measurements.
For a linear subspace $L\subset{\mathbb R}^V$, let $L^\perp$ be the
orthogonal complement of $L$ and let $P_L$ be the
orthogonal projector onto the subspace $L$. Denote $L:=\operatorname{
supp}(S_{\ast})$, $r=\operatorname{ rank}(S_{\ast})$. A \emph{coherence
coefficient} is a constant $\nu\geq1$ such that
%
\begin{eqnarray}
\label{coherA} \qquad\|P_L e_v\|^2 &\leq&
\frac{\nu r }{m},\qquad v\in V \quad\mbox{and}
\nonumber
\\[-8pt]
\\[-8pt]
\nonumber
\bigl |\bigl\langle\operatorname{
sign}(S_{\ast})e_u, e_v\bigr
\rangle\bigr|^2 &\leq&\frac{\nu
r}{m^2}, \qquad u,v\in V
\end{eqnarray}
(it is easy to see that $\nu$ cannot be smaller than $1$).

The following highly nontrivial result is essentially due to Candes and Tao
\cite{CandesTao}
(a version stated here is due to Gross~\cite{Gross-2} and it is an
improvement of the initial result of Candes and Tao).
It shows that target matrices of ``low coherence''
(for which $\nu$ is a relatively small constant) can be recovered
exactly using
the nuclear norm minimization algorithm (\ref{nucemin}) provided that
the number of observed entries is of the order $mr$ (up to a log factor).

\begin{theorem}
\label{gross}
Suppose conditions (\ref{coherA}) hold for some
$\nu\geq1$. Then, there exists a numerical constant $C>0$ such that,
for all
$n\geq C\nu r m \log^2 m$, $\hat S=S_{\ast}$ with probability at
least $1-m^{-2}$.
\end{theorem}

In the case of noisy matrix completion, a matrix version of LASSO is
based on a trade-off between fitting the target matrix to the data
using least squares and minimizing the nuclear norm
%
\begin{equation}
\label{LSnuce} \hat S:=\mathop{\operatorname{argmin}}_{S\in{\cal S}_V}
\Biggl[n^{-1}\sum_{j=1}^n
\bigl(Y_j-S\bigl(X_j,X_j'
\bigr)\bigr)^2 +\eps\|S\|_1 \Biggr].
\end{equation}
This method and its modifications have been studied by a number of authors;
see~\cite{CandesPlan,Rohde,Negahban,Ko26,Ko27}. The
following low-rank oracle inequality was proved in~\cite{Ko26}
(Theorem 4) for a ``linearized version'' of the matrix LASSO estimator~$\hat S$.
Assume that, for some constant $a>0$, $|Y|\leq a$ a.s.
Let $t>0$ and suppose that
$
\eps\geq4 a (\sqrt{\frac{t+\log(2m)}{nm}}\vee\frac
{2(t+\log(2m))}{n} ).
$
Then, there exists a constant $C>0$ such that with probability at least
$1-e^{-t}$
\[
\|\hat S-S_{\ast}\|_{L_2(\Pi^2)}^2 \leq \inf
_{S\in{\cal S}_V} \bigl[\|S-S_{\ast}\|_{L_2(\Pi^2)}^2+C
m^2 \eps^2 \operatorname{ rank}(S) \bigr].
\]
In particular,
$
\|\hat S-S_{\ast}\|_{L_2(\Pi^2)}^2 \leq
C m^2 \eps^2 \operatorname{ rank}(S_{\ast}).
$
Very recently, the last bound was proved in~\cite{Klopp}
for the matrix LASSO estimator (\ref{LSnuce}) itself in the case when
the domain
of optimization problem is $\{S\dvtx \|S\|_{L_{\infty}}\leq a\}$, where
$\|S\|_{L_{\infty}}:=\max_{u,v\in V}|S(u,v)|$; in fact, both~\cite{Ko26}
and~\cite{Klopp} dealt with the case of rectangular matrices.

In the current paper, we are more interested in
the case when the target kernel~$S_{\ast}$ is defined on the set $V$
of vertices
of a weighted graph $G=(V,A)$ with a symmetric matrix
$A:=(a(u,v))_{u,v\in V}$ of nonnegative weights.
This allows one to define
the notion of graph Laplacian and to introduce discrete Sobolev norms
characterizing smoothness of functions on $V$ as well as symmetric kernels
on $V\times V$.
Denote
$
\operatorname{ deg}(u):=\sum_{v\in V}a(u,v), u\in V.
$
It is common in graph theory to call $\operatorname{ deg}(u)$ the degree of
vertex $u$.
Let $D$ be the diagonal $m\times m$ matrix (kernel) with the degrees of vertices
on the diagonal (it is assumed that the vertices of the graph have been ordered
in an arbitrary, but fixed way). The Laplacian of the weighted graph
$G$ is
defined as $\Delta:=D-A$. Denote $\langle\cdot,\cdot\rangle$ the canonical
Euclidean inner product in the $m$-dimensional space ${\mathbb R}^V$ of
functions $f\dvtx V\mapsto{\mathbb R}$ and let $\|\cdot\|$ be the
corresponding norm. It is easy to see that
\[
\langle\Delta f,f\rangle= \frac{1}{2}\sum_{u,v\in V}
a(u,v) \bigl(f(u)-f(v)\bigr)^2,
\]
implying that $\Delta\dvtx  {\mathbb R}^{V}\mapsto{\mathbb R}^V$ is a
symmetric nonnegatively definite linear transformation. In a special
case of a usual
graph $(V,E)$ with vertex set $V$ and edge set $E$, one defines $A(u,v)=1$
if and only if $u\sim v$ (i.e., vertices $u$ and $v$ are connected with
an edge)
and $A(u,v)=0$ otherwise. In this case, $\operatorname{ deg}(u)$ is the number
of edges incident to the vertex $u$ and
$
\langle\Delta f,f\rangle= \sum_{u\sim v} (f(u)-f(v))^2.
$
The notion of graph Laplacian allows one to define discrete Sobolev
norms $\|\Delta^{q/2}f\|, q>0$
for functions on the vertex set of the graph and thus to describe their
smoothness on the graph. Given a symmetric kernel $S\dvtx V\times V\mapsto
{\mathbb R}$, one can also describe its smoothness in terms of the
norms $\|\Delta^{q/2}S\|_2$.
Suppose $S$ has the following spectral representation:
$
S=\sum_{j=1}^m \mu_j (\psi_j \otimes\psi_j),
$
where $\mu_j, j=1,\ldots, m$ are the eigenvalues of $S$ (repeated with their
multiplicities) and $\psi_j, j=1,\ldots, m$ are the corresponding orthonormal
eigenfunctions in ${\mathbb R}^V$, then
\begin{eqnarray*}
\bigl\|\Delta^{q/2}S\bigr\|_2^2 &=& \operatorname{ tr}
\bigl(\Delta^{q/2}S^2 \Delta^{q/2}\bigr)=
\operatorname{ tr}\bigl(\Delta^q S^2\bigr)=\sum
_{j=1}^m \mu_j^2\bigl
\langle\Delta^q \psi_j, \psi_j\bigr\rangle\\
&=&
\sum_{j=1}^m \mu_j^2
\bigl\|\Delta^{q/2}\psi_j\bigr\|^2.
\end{eqnarray*}
Basically, it means that the smoothness of the kernel $S$ depends on the
smoothness of its eigenfunctions. In what follows, we will often use rescaled
versions of Sobolev norms,
\[
\bigl\|\Delta^{q/2}f\bigl\|_{L_2(\Pi)}=m^{-1/2}\bigl\|
\Delta^{q/2}f\bigr\|^2,\qquad \bigl\|\Delta^{q/2}S
\bigr\|_{L_2(\Pi^2)}=m^{-1}\bigl\|\Delta^{q/2}S\bigr\|_2.
\]

It will be convenient for our purposes to fix $q>0$
and to define a nonnegatively definite symmetric kernel $W:=\Delta^q$.
We will characterize the smoothness of a kernel
$S\in{\cal S}_V$ by the squared Sobolev-type norm $\|W^{1/2}S\|
_{L_2(\Pi^2)}^2$.
The kernel $W$ will be fixed throughout the paper, and its spectral
properties are crucial
in our analysis.\footnote{In fact, the relationship of $W$ to the
graph and its Laplacian
will be of little importance allowing, possibly, other interpretations
of the problem.}
Assume that $W$ has the following spectral representation
$
W=\sum_{k=1}^m \lambda_k (\phi_k \otimes\phi_k),
$
where $0\leq\lambda_1 \leq\cdots\leq\lambda_m$ are the eigenvalues
repeated with their multiplicities, and $\phi_1,\ldots, \phi_m$ are
the corresponding orthonormal eigenfunctions (of course, there is a
multiple choice of $\phi_k$ in the case of repeated eigenvalues). Let
$k_0:=\min\{k\leq m\dvtx \lambda_k>0\}$.
We will assume in what follows that, for some constant
$c\geq1$, $\lambda_{k+1}\leq c\lambda_k$ for all $k\geq k_0$.
It will be also convenient to set $\lambda_k:=+\infty, k>m$.

Let $\rho:=\|W^{1/2}S_{\ast}\|_{L_2(\Pi^2)}$ and $r:=\operatorname{
rank}(S_{\ast})$.
It is easy to show
(see the proof of Theorem~\ref{LSupth} below) that kernel $S_{\ast}$
can be approximated by the following kernel:
$
S_{\ast,l}:=\sum_{i,j=1}^{l} \langle S_{\ast} \phi_i, \phi
_j\rangle(\phi_i\otimes\phi_j)
$
with the approximation error
%
\begin{equation}
\label{apprerr} \|S_{\ast}-S_{\ast,l}\|_{L_2(\Pi^2)}^2
\leq\frac{2\rho^2}{\lambda_{l+1}}.
\end{equation}
Note that the kernel $S_{\ast,l}$
can be viewed as an $l\times l$ matrix (represented in the basis of
eigenfunctions $\{\phi_j\}$)
and $\operatorname{ rank} (S_{\ast,l})\leq r\wedge l$,
so, one needs $\sim(r\wedge l)l$ parameters to characterize such matrices.
Thus, one can expect, that such a kernel can be estimated, based on $n$ linear
measurements, with the squared $L_2(\Pi^2)$-error of the order
$\frac{a^2(r\wedge l)l}{n}$. Taking into account the bound on the
approximation error (\ref{apprerr}) and optimizing
with respect to $l=1,\ldots, m$, it would be also natural to expect the
following error rate in the problem of estimation of the target kernel
$S_{\ast}\dvtx $
%
\begin{equation}
\label{upperbd} \min_{1\leq l\leq m} \biggl[\frac{a^2(r\wedge l)l}{n}\vee
\frac
{\rho^2}{\lambda_{l+1}} \biggr].
\end{equation}
We will show that such a rate is attained (up to constants and log factors)
for a version of least squares method with a nonconvex complexity
penalty; see Section~\ref{sectLSnonconvex}.
This method is not computationally tractable, so, we also study another
method, based on convex penalization with a combination of nuclear norm
and squared Sobolev type norm, and show that the rates are attained
for such a method, too, provided that the target matrix satisfies a version
low coherence assumption with respect to the basis of eigenfunctions
of $W$. More precisely, we will prove error bounds involving
so called \emph{coherence function}
$\varphi(S_{\ast};\lambda):= \langle P_{\operatorname{ supp}(S_{\ast})},
\sum_{\lambda_j\leq\lambda}(\phi_j\otimes\phi_j) \rangle$,
that characterizes the relationship between the kernel $W$ defining the
smoothness and the target kernel
$S_{\ast}$; see Section~\ref{sectLSconvex} for more details; see
also~\cite{KoRan} for similar results in the case of ``linearized
least squares'' estimator with double penalization.
Finally, we prove minimax lower
bounds on the error rate that are roughly of the order
$
\max_{1\leq l\leq m} [\frac{a^2(r\wedge l)l}{n}\wedge\frac
{\rho^2}{\lambda_{l}} ]
$
(subject to some extra conditions and with additional terms; see
Section~\ref{sectlowerbound}).
In typical situations, this expression is, up to a constant, of the
same order as the upper bound
(\ref{upperbd}).
For instance, if $\lambda_l\asymp l^{2\beta}$ for some $\beta>1/2$,
then the minimax error rate of estimation of the target kernel $S_{\ast}$
is of the order
\[
\biggl( \biggl(\frac{a^2 \rho^{1/\beta} r}{n} \biggr)^{2\beta/(2\beta+1)} \wedge
\biggl(\frac{a^2\rho^{2/\beta}}{n} \biggr)^{\beta/(\beta+1)} \wedge
\frac{a^2 rm}{n} \biggr)\vee\frac{a^2}{n}
\]
(up to log factors). When $m$ is sufficiently large, the term $\frac
{a^2 rm}{n}$ will be dropped from the minimum, and we end up with a
nonparametric
convergence rate controlled by the smoothness parameter $\beta$ and
the rank
$r$ of the target matrix $S_{\ast}$ (the dependence on $m$ in the
first two
terms of the minimum is only in the log factors).

The focus of the paper is on the matrix completion problems with uniform
random design, but it is very straightforward to extend the results of the
following sections to sampling models with more general design
distributions discussed in the literature on low rank recovery (such
as, e.g., the models of random linear\vadjust{\goodbreak} measurements studied in
\cite{Ko26,Ko27}). It is also not hard to replace the range $a$ of
the response
variable $Y$ by the standard deviation of the noise in the upper and
lower bounds
obtained below. This is often done in the literature on low-rank
recovery, and
it can be easily extended to the framework discussed in the paper by modifying
our proofs. We have not discussed this in the paper due to the lack of space.

\section{Minimax lower bounds}\label{sectlowerbound}

In this section, we derive minimax lower bounds on the $L_2(\Pi
^2)$-error of an arbitrary estimator $\hat S$ of the target kernel
$S_{\ast}$ under the assumptions
that the response variable $Y$ is bounded by a constant $a>0$, the rank
of $S_{\ast}$ is bounded by $r\leq m$ and its Sobolev norm $\|
W^{1/2}S_{\ast}\|_{L_2(\Pi^2)}$ is
bounded by $\rho>0$. More precisely, given $r=1,\ldots, m$ and $\rho
>0$, denote\vspace*{1pt}
by ${\cal S}_{r,\rho}$ the set of all symmetric kernels
$S\dvtx V\times V\mapsto{\mathbb R}$ such that $\operatorname{ rank}(S)\leq r$
and \mbox{$\|W^{1/2}S\|_{L_2(\Pi^2)}\leq\rho$}.
Given $r, \rho$ and $a>0$, let ${\cal P}_{r,\rho,a}$ be the set of
all probability distributions of $(X,X',Y)$ such that
$(X,X')$ is uniformly distributed in $V\times V$,
$|Y|\leq a$ a.s. and
${\mathbb E}(Y|X,X')=S_{\ast}(X,X')$,
where $S_{\ast}\in{\cal S}_{r,\rho}$.
For $P\in{\cal P}_{r,\rho, a}$,
denote $S_P(X,X'):={\mathbb E}_P (Y|X,X')$.

Recall that $\{\phi_j, j=1,\ldots, m\}$ are the eigenfunctions of $W$
orthonormal in the space $({\mathbb R}^V, \langle\cdot,\cdot\rangle)$.
Then $\bar\phi_j:=\sqrt{m} \phi_j, j=1,\ldots, m$ are orthonormal
in~$L_2(\Pi)$.

We will obtain minimax lower bounds for classes of distributions ${\cal
P}_{r,\rho,a}$ in two different cases. Define
$
Q_p:=\max_{1\leq j\leq m}\|\bar\phi_j\|_{L_p(\Pi)}^2.
$
In the first case, we assume that,
for some (relatively large) value of $p\geq2$, the quantity $Q_p$ is not
too large. Roughly, it means that most of the components of vectors
$\phi_j\in{\mathbb R}^V$ are uniformly small, say,
$
\phi_j(v)\asymp m^{-1/2}, v\in V, j=1,\ldots, m.
$
In other words, the $m\times m$ matrix $(\phi_j(v))_{j=1,\ldots,
m;v\in V}$
is ``dense,'' so we refer to this case as a \textit{``dense case.''}
The opposite case is when this matrix is \textit{``sparse.''} Its ``sparsity''
will be characterized by the quantity
\[
d:=\max_{v\in V}\operatorname{ card}\bigl\{j\dvtx
\phi_j(v)\neq0\bigr\},
\]
which, in this case, should be relatively small.
A typical example is the case when basis of eigenfunctions $\{\phi
_j,j=1,\ldots, m\}$ coincides with the canonical basis $\{e_v\dvtx v\in V\}$
of ${\mathbb R}^V$ (then, $d=1$).

Denote $l_0:=k_0\wedge32$.
In the \textit{dense case}, the following theorem holds.

\begin{theorem}
\label{lowerbound}
Define
\[
\delta_n^{(1)}(r,\rho,a):= \max_{l_0\leq l\leq m}
\biggl[\frac{a^2(r\wedge l)l}{n} \wedge \frac{\rho^2}{\lambda_l} \wedge
\frac{1}{p-1}\frac{1}{Q_p^2}\frac{a^2(r\wedge l)}{l}\frac
{1}{m^{4/p}} \biggr].
\]
There exist constants $c_1, c_2>0$ such that
\[
\inf_{\hat S_n}\sup_{P\in{\cal P}_{r,\rho,a}}{\mathbb
P}_{P} \bigl\{ \|\hat S_n-S_P
\|_{L_2(\Pi^2)}^2\geq c_1 \delta_n^{(1)}(r,
\rho, a) \bigr\}\geq c_2,
\]
where the infimum is taken over all the estimators $\hat S_n$ based on
$n$ i.i.d.
copies of $(X,X',Y)$.
\end{theorem}

In fact, it will follow from the proof that, if $\lambda_{k_0}\leq
\frac{n\rho^2}{a^2(r\wedge k_0)k_0}$ (i.e., the smallest nonzero
eigenvalue of $W$ is not too large), then the maximum in the definition
of $\delta_n^{(1)}(r,\rho,a)$ can be extended to all $l=1,\ldots, m$.

\begin{corollary}
\label{lowerboundcorr}
Let
\[
\delta_n^{(2)}(r,\rho,a):= \max_{l_0\leq l\leq m}
\biggl[\frac{a^2(r\wedge l)l}{n} \wedge \frac{\rho^2}{\lambda_l} \wedge
\frac{1}{Q_{\log m}^2}\frac{a^2(r\wedge l)}{l}\frac{1}{\log m} \biggr].
\]
There exist constants $c_1, c_2>0$ such that
\[
\inf_{\hat S_n}\sup_{P\in{\cal P}_{r,\rho,a}}{\mathbb
P}_{P} \bigl\{ \|\hat S_n-S_P
\|_{L_2(\Pi^2)}^2\geq c_1 \delta_n^{(2)}(r,
\rho, a) \bigr\}\geq c_2.
\]
\end{corollary}

\begin{pf} Take $p=\log m$ in the statement of Theorem~\ref{lowerbound}
and observe that $m^{4/p}=e^{4}$ and $\frac{1}{p-1}\geq\frac{1}{\log m}$.
\end{pf}

\begin{remark*} It is easy to check that $e^{-2}Q_{\infty}\leq Q_{\log
m}\leq Q_{\infty}$.
\end{remark*}

It is obvious that one can replace
the quantity $\delta_n^{(1)}(r,\rho, a)$ in Theorem~\ref{lowerbound}
[or the quantity $\delta_n^{(2)}(r,\rho, a)$ in Corollary \ref
{lowerboundcorr}]
by the following smaller quantity:
\[
\delta_n^{(3)}(r,\rho,a):= \max_{l_0\leq l\leq L}
\biggl[\frac{a^2(r\wedge l)l}{n} \wedge \frac{\rho^2}{\lambda_l} \biggr],
\]
where
$
L:= [\frac{1}{Q_p m^{2/p}}\sqrt{\frac{n}{p-1}} ]\wedge m.
$
Moreover, denote
\[
\bar l:= \max \biggl\{l=l_0,\ldots, m\dvtx (r\vee l)l
\lambda_l \leq \frac{\rho^2 n}{a^2} \biggr\}.
\]
It is straightforward to check that
\[
\max_{l_0\leq l\leq m} \biggl[\frac{a^2(r\wedge l)l}{n} \wedge
\frac{\rho^2}{\lambda_l} \biggr] = \frac{a^2(r\wedge\bar l)\bar l}{n} \vee\frac{\rho
^2}{\lambda_{\bar l+1}}
\]
and, if $\bar l\leq L$, then
$
\delta_n^{(3)}(r,\rho,a)=\frac{a^2(r\wedge\bar l)\bar l}{n} \vee
\frac{\rho^2}{\lambda_{\bar l+1}}.
$

\begin{example*} Suppose that, for some $\beta>1/2$, $\lambda_l\asymp
l^{2\beta},
l=1, \ldots, m$ (in particular, it means that $\lambda_l\neq0$ and
$l_0=k_0=1$).
Then, an easy computation shows that
\[
\bar l = (\check l \wedge m)\vee1,\qquad \check l \asymp \biggl(\frac{\rho^2}{a^2}
\frac{n}{r} \biggr)^{1/(2\beta+1)}\wedge \biggl(
\frac{\rho^2 n}{a^2} \biggr)^{1/(2\beta+2)}.
\]
Let $p=\log m$ and take
$
L:= [\frac{1}{e^2 Q_p}\sqrt{\frac{n}{\log(m/e)}}
]\wedge m.
$
The condition $\bar l\leq L$ is satisfied, for instance, when either
$e^2 Q_p \sqrt{\log(m/e)}  (\frac{\rho^2}{a^2 r}
)^{1/(2\beta+1)}
\leq c' n^{{1}/{2}-{1}/{(2\beta+1)}}$,
or
$e^2 Q_p \sqrt{\log(m/e)}  (\frac{\rho}{a} )^{1/(\beta+1)}
\leq c' n^{{1}/{2}-{1}/{(2\beta+2)}}$,
where $c'>0$ is a small enough constant
(this, essentially, means that
$n$ is sufficiently large). Under either of these conditions, we get
the following expression for a minimax lower bound:
%
\begin{equation}
\label{minmaxlower} \biggl( \biggl(\frac{a^2 \rho^{1/\beta} r}{n} \biggr)^{2\beta/(2\beta+1)}
\wedge \biggl(\frac{a^2 \rho^{2/\beta}}{n} \biggr)^{\beta
/(\beta+1)} \wedge
\frac{a^2 rm}{n} \biggr) \vee\frac{a^2}{n}.
\end{equation}
\end{example*}

We now turn to the \textit{sparse case}.

\begin{theorem}
\label{lowerboundsparse}
Let
\[
\delta_n^{(4)}(r,\rho,a):= \max_{l_0\leq l\leq m}
\biggl[\frac{a^2(r\wedge l)l}{n} \wedge \frac{\rho^2}{\lambda_l} \wedge
\frac{a^2}{d \log m}\frac{l^2}{m^2} \biggr].
\]
There exist constants $c_1, c_2>0$ such that
\[
\inf_{\hat S_n}\sup_{P\in{\cal P}_{r,\rho,a}}{\mathbb
P}_{P} \bigl\{ \|\hat S_n-S_P
\|_{L_2(\Pi^2)}^2\geq c_1 \delta_n^{(4)}(r,
\rho, a) \bigr\}\geq c_2.
\]
\end{theorem}

It will be clear from the upper bounds of Section~\ref{sectLSnonconvex}
(see the remark after Theorem~\ref{LSupth}) that,
at least in a special case when $\{\phi_j\}$ coincides with the
canonical basis
of ${\mathbb R}^V$, the additional term $\frac{a^2}{d \log m}\frac
{l^2}{m^2}$ is correct (up to a log factor). At the same time,
most likely, the ``third terms'' of the bounds of Theorem~\ref{lowerbound}
(in the dense case) and Theorem~\ref{lowerboundsparse} (in the
sparse case)
have not reached their final form yet. A more sophisticated construction
of ``well separated'' subsets of ${\cal P}_{r,\rho,a}$ might be needed
to achieve this goal. The main difficulty in the proof given below is
related to the fact that we have to impose constraints, on the one hand,
on the entries of the target matrix represented in the canonical basis
and, on the other hand, on the Soblolev type norm $\|W^{1/2}S\|
_{L_2(\Pi^2)}$
(for which it is convenient to use the representation in the basis of
eigenfunctions of $W$). Due to this fact, we are using the last representation
in our construction, and we have to use an argument based
on the properties of Rademacher sums to ensure that the entries of the matrix
represented in the canonical basis are uniformly bounded by $a$. This
is the reason why the ``third terms'' occur in the bounds of Theorems
\ref{lowerbound} and~\ref{lowerboundsparse}. In this case, when the
constraints are only on the norm $\|W^{1/2}S\|_{L_2(\Pi^2)}$ and on the
variance of the noise and there are no constraints on $\|S\|_{L_{\infty}}$,
it is much easier to prove the lower bound of the order
$
\max_{l_0\leq l\leq m} [\frac{\sigma^2(r\wedge l)l}{n}
\wedge \frac{\rho^2}{\lambda_l}
]
$
without any additional terms. Note, however, that the condition
$\|S_{\ast}\|_{L_{\infty}}\leq a$ is of importance in the following sections
to obtain the upper bounds for penalized least squares estimators that match
the lower bounds up to log factors.

\begin{pf*}{Proof of Theorem~\ref{lowerbound}} The proof relies on several
well-known facts stated below.
In what follows, $K(\mu\|\nu):=-{\mathbb E}_{\mu}\log\frac{d\nu
}{d\mu}$ denotes Kullback--Leibler divergence
between two probability measures $\mu,\nu$ defined on the same space
and such that $\nu\ll\mu$ (i.e., $\nu$ is absolutely continuous
with respect to $\mu$). We will denote by $P^{\otimes n}$ the $n$-fold
product measure $P^{\otimes n}:=P\otimes P\cdots\otimes P$.
The following proposition is a version of Theorem 2.5 in~\cite{Tsybakov}.

\begin{proposition}\label{KBound}
Let ${\cal P}$ be a finite set of distributions
of $(X,X',Y)$ such that the following assumptions hold:
\begin{longlist}[(1)]
\item[(1)] there exists
$P_0\in{\cal P}$ such that
for all $P \in{\cal P}$, $P \ll P_0$;
\item[(2)] there exists $\alpha\in(0,1/8)$ such that
\[
\sum_{P \in{\cal P}} K\bigl(P_0^{\otimes n}
\|P^{\otimes n}\bigr) \leq\alpha \bigl(\operatorname{ card}({\cal P})-1\bigr)
\log\bigl( \operatorname{ card}({\cal P})-1\bigr);
\]
\item[(3)] for all $P_1, P_2 \in{\cal P}$, $\|S_{P_1} - S_{P_2}\|
_{L_2(\Pi^2)}^2 \geq4s^2 > 0$.
\end{longlist}
Then, there exists a constant $\beta>0$ such that
%
\begin{equation}
\inf_{\hat{S}_n} \max_{P \in{\cal P}} {\mathbb
P}_P\bigl\{\|\hat {S}_n-S_P
\|_{L_2(\Pi^2)}^2 \geq s^2\bigr\} \geq\beta>0.
\end{equation}
\end{proposition}

We will also use Varshamov--Gilbert bound (see~\cite{Tsybakov}, Lemma 2.9,
page~104), Sauer's lemma (see~\cite{Ko27}, page 39)
and the following elementary bound for Rademacher sums
(\cite{Gine}, page 21):
for all $p\geq2$,
%
\begin{equation}
\label{Rademach} {\mathbb E}^{1/p} \Biggl|\sum_{j=1}^N
\eps_j t_j \Biggr|^p \leq\sqrt{p-1} \Biggl(\sum
_{j=1}^N t_j^2
\Biggr)^{1/2},\qquad (t_1,\ldots, t_N)\in{\mathbb
R}^N,
\end{equation}
where
$\eps_1, \ldots, \eps_N$ are i.i.d. Rademacher random variables
(i.e., $\eps_j=+1$ with probability $1/2$ and $\eps_j=-1$ with
the same probability).

We will start the proof with constructing a ``well separated'' subset
${\cal P}$
of the class of distributions ${\cal P}_{r,\rho, a}$ that will allow
us to
use Proposition~\ref{KBound}. Fix $l\leq m$, $l\geq32$ and $\kappa
>0$. Denote $l'=[l/2], l''=l-l'$.
First assume that $r\leq l''$. Denote
$R_{\sigma}:=\kappa ((\sigma_{ij})\dvtx i=1,\ldots, l', j=1,\ldots,
r )$, where $\sigma_{ij}=+1$ or $\sigma_{ij}=-1$.
Let ${\cal R}_{l',r}=\{R_{\sigma}\dvtx \sigma\in\{-1,1\}^{l'\times r}\}$
(so, ${\cal R}_{l',r}$ is the class of all $l'\times r$ matrices with
entries $+\kappa$ or $-\kappa$).
Given
$R\in{\cal R}_{l',r}$, let
\[
\tilde R:= \pmatrix{
 R & R & \cdots& R &
O_{l',l^*}}
\]
be the $l'\times l''$ matrix that consists of $[l''/r]$ blocks $R$
and the last block $O_{l',l^*}$, where
$l^{\ast}:=l''-[l''/r]r$ and $O_{k_1,k_2}$ is
the $k_1\times k_2$ zero matrix.
Finally, define the following symmetric $m\times m$ matrix:
\[
R^{\diamondsuit}:=\pmatrix{ O_{l',l'}
& \tilde R & O_{l',m-l}
\vspace*{2pt}\cr
\tilde R^{T} & O_{l'',l''} & O_{l'',m-l}
\vspace*{2pt}\cr
O_{m-l,l'} & O_{m-l,l''} & O_{m-l,m-l}}.
\]
Now, given $\sigma\in\{-1,1\}^{l'\times r}$, define a
symmetric kernel $K_{\sigma}\dvtx V\times V\mapsto{\mathbb R}$,
\[
K_{\sigma}:=\sum_{i,j=1}^m
\bigl(R_{\sigma}^{\diamondsuit}\bigr)_{ij} (\phi _i
\otimes\phi_j).
\]
It is easy to see that
%
\begin{eqnarray}
\label{Ksigma}  K_{\sigma}(u,v)&=&K_{\sigma}^{\prime }(u,v)+K_{\sigma}^{\prime }(v,u),
\nonumber
\\[-8pt]
\\[-8pt]
\nonumber
K_{\sigma}^{\prime }(u,v)&=& \kappa\sum
_{i=1}^{l'}\sum_{j=1}^r
\sigma_{ij}\phi_i(u)\sum_{k=0}^{[l''/r]-1}
\phi_{l'+rk+j}(v).
\end{eqnarray}

Let
$
\Lambda:= \{\sigma\in\{-1,1\}^{l'\times r}\dvtx \max_{u,v\in
V}|K_{\sigma}(u,v)|\leq a\}.
$
We will show that, if $\kappa$ is sufficiently small (its precise value
to be specified later), then
the set $\Lambda$ contains at least three quarters of the points
of the combinatorial cube $\{-1,1\}^{l'\times r}$.
To this end, define
$
\xi:= \max_{u,v\in V}|K_{\eps}(u,v)|,
$
where $\eps\in\{-1,1\}^{l'\times r}$ is a random vector with
i.i.d. Rademacher components. Assume, in addition, that
$\eps$ and $(X,X')$ are independent.
It is enough to show that $\xi\leq a$ with probability
at least $3/4$.
We have
%
\begin{eqnarray}
\label{ogogo}  {\mathbb P}\{\xi\geq a\}&\leq&\sum_{u,v\in V}{
\mathbb P}\bigl\{\bigl|K_{\eps
}(u,v)\bigr|\geq a\bigr\}
\nonumber\\
&
=&m^2 {\mathbb E} {\mathbb P}\bigl\{\bigl|K_{\eps}
\bigl(X,X'\bigr)\bigr|\geq a|X,X'\bigr\}
\\
&=&m^2 {\mathbb P}\bigl\{\bigl|K_{\eps}
\bigl(X,X'\bigr)\bigr|\geq a\bigr\}\leq \frac{m^{2}{\mathbb
E}|K_{\eps}(X,X')|^p}{a^p}.\nonumber
\end{eqnarray}
We will use bound (\ref{Rademach}) to control ${\mathbb E}(|K_{\eps
}(X,X')|^p|X,X')$
[recall that $K_{\eps}(u,v),\break  u,v\in V$ is a Rademacher sum].
Denote
\[
\tau^2(u,v):= \sum_{i=1}^{l'}
\sum_{j=1}^r \phi_i^2(u)
\Biggl(\sum_{k=0}^{[l''/r]-1}\phi_{l'+rk+j}(v)
\Biggr)^2. 
\]
Observe that
$
\tau^2(u,v)\leq\frac{l''}{r}q(l',u)q(l'',v)\leq
q(l,u)q(l,v)\frac{l}{r},
$
where
$
q(l,u):=\break\sum_{j=1}^l \phi_j^2(u), u\in V,
$
and we used the bound
%
\begin{equation}
\label{lr-d} \Biggl(\sum_{k=0}^{[l''/r]-1}
\phi_{l'+rk+j}(v) \Biggr)^2\leq \frac{l''}{r} \sum
_{k=0}^{[l''/r]-1}\phi_{l'+rk+j}^2(v).
\end{equation}
Thus, applying (\ref{Rademach}) to the Rademacher sum $K_{\eps}^{\prime }$,
we get
\begin{eqnarray*}
{\mathbb E}\bigl|K_{\eps}(u,v)\bigr|^p&\leq&2^{p-1} \bigl({
\mathbb E}\bigl|K_{\eps
}^{\prime }(u,v)\bigr|^p+{\mathbb
E}\bigl|K_{\eps}^{\prime }(v,u)\bigr|^p \bigr)\\
& \leq&
2^p (p-1)^{p/2} \kappa^p \bigl(
\tau^2(u,v)\vee\tau^2(v,u)\bigr)^{p/2}\\
&\leq&
2^p (p-1)^{p/2} \kappa^p q^{p/2}(l,u)q^{p/2}(l,v)
\biggl(\frac
{l}{r} \biggr)^{p/2}.
\end{eqnarray*}
Given $p\in[2,+\infty]$, denote
$
Q_p(l):=  \|\frac{m}{l}q(l,\cdot) \|_{L_{p/2}(\Pi)}=
\|\frac{1}{l}\sum_{j=1}^l \bar\phi_j^2 \|_{L_{p/2}(\Pi)}$
for $l=1,\ldots, m$.
This yields
\begin{eqnarray*}
{\mathbb E}\bigl|K_{\eps}\bigl(X,X'
\bigr)\bigr|^p&=& {\mathbb E} {\mathbb E}\bigl(\bigl|K_{\eps}
\bigl(X,X'\bigr)\bigr|^p|X,X'\bigr)\\
&\leq&
2^p (p-1)^{p/2}\kappa^p \biggl(
\frac{l}{r} \biggr)^{p/2}{\mathbb E}\bigl(q^{p/2}(l,X)q^{p/2}
\bigl(l,X'\bigr)\bigr)\\
&=&
2^p (p-1)^{p/2}\kappa^p \biggl(
\frac{l}{r} \biggr)^{p/2} \bigl({\mathbb E}q^{p/2}(l,X)
\bigr)^2\\
&=&
2^p (p-1)^{p/2}\kappa^p \biggl(
\frac{l}{r} \biggr)^{p/2} \biggl(\frac{l}{m}
\biggr)^p Q_p^{p}(l).
\end{eqnarray*}
Substituting the last bound into (\ref{ogogo}), we get
\[
{\mathbb P}\{\xi\geq a\}\leq \frac{m^{2}{\mathbb E}|K_{\eps}(X,X')|^p}{a^p} \leq
m^2 2^p(p-1)^{p/2}\frac{\kappa^p}{a^p} \biggl(
\frac{l}{r} \biggr)^{p/2} \biggl(\frac{l}{m}
\biggr)^p Q_p^{p}(l).
\]
Now, to get ${\mathbb P}\{\xi\geq a\}\leq1/4$, it is enough to take
%
\begin{equation}
\label{kappa1} \kappa\leq 2^{-(1+2/p)}(p-1)^{-1/2}
\frac{1}{Q_p(l)}\frac{m}{l}\frac{a\sqrt {r}}{\sqrt{l}} \frac{1}{m^{2/p}}.
\end{equation}

Next observe that
\[
\operatorname{ card}(\Lambda)\geq\frac{3}{4} 2^{l'r} >
\sum_{k=0}^{[l'r/2]}\pmatrix{l'r
\cr k}.
\]
It follows from Sauer's lemma that there exists a subset $J\subset\{
(i,j)\dvtx 1\leq i\leq l', 1\leq j\leq r\}$ with $\operatorname{ card}(J)=[l'r/2]+1$
and such that
$\pi_J(\Lambda)=\{-1,1\}^{J}$, where $\pi_J\dvtx \{-1,1\}^{l'\times
r}\mapsto\{-1,1\}^J$,
$
\pi_J  (\sigma_{ij}\dvtx i=1,\ldots, l', j=1,\ldots, r )=
(\sigma_{ij}\dvtx\break (i,j)\in J ).
$
Since $l\geq32$, we have $l'r\geq16$ and $\operatorname{ card}(J)\geq8$. We
can now apply Varshamov--Gilbert bound to the combinatorial cube
$\{-1,1\}^J$ to prove that there exists a subset $E\subset\{-1,1\}^J$
such that $\operatorname{ card}(E)\geq2^{l'r/16}+1$ and, for all $\sigma',
\sigma''\in E,
\sigma'\neq\sigma''$,
$
\sum_{(i,j)\in J}I(\sigma^{\prime }_{ij}\neq\sigma^{\prime \prime}_{ij})\geq\frac{l'r}{16}.
$
It is now possible to choose a subset $\Lambda'$ of $\Lambda$ such that
$\operatorname{ card}(\Lambda')=\operatorname{ card}(E)$ and $\pi_J (\Lambda')=E$.
Then, we have
$\operatorname{ card}(\Lambda')\geq2^{l'r/16}+1$ and
%
\begin{equation}
\label{hemming} \sum_{i=1}^{l'}\sum
_{j=1}^r I\bigl(\sigma^{\prime }_{ij}
\neq\sigma ^{\prime \prime}_{ij}\bigr)\geq\frac{l'r}{16}
\end{equation}
for all $\sigma', \sigma''\in\Lambda', \sigma'\neq\sigma''$.

We are now in a position to define the set of distributions ${\cal P}$.
For $\sigma\in\Lambda'$, denote by $P_{\sigma}$ the distribution
of $(X,X',Y)$ such that $(X,X')$ is uniform in $V\times V$ and\vadjust{\goodbreak}
the conditional distribution of $Y$ given $(X,X')$ is defined as follows:
\[
{\mathbb P}_{P_{\sigma}}\bigl\{Y= \delta a|X,X'\bigr\}
=p_{\sigma}\bigl(X,X'\bigr)=1/2+ \delta K_{\sigma}
\bigl(X,X'\bigr)/8a,\qquad \delta\in\{ -1,+1\}.
\]

Since $|K_{\sigma}(X,X')|\leq a$ for all $\sigma\in\Lambda'$,
we have $p_{\sigma}(X,X')\in[3/8,5/8],\sigma\in\Lambda$.
Denote ${\cal P}:=\{P_{\sigma}\dvtx \sigma\in\Lambda'\}$.
For $P=P_{\sigma}\in{\cal P}$, we have
\[
S_P(u,v)={\mathbb E}\bigl(Y|X=u,X'=v\bigr)=
\tfrac{1}{4}K_{\sigma}(u,v).
\]
Note that
$
\operatorname{ rank}(S_P)=\operatorname{ rank}(K_{\sigma})=\operatorname{ rank}(R_{\sigma
}^{\diamondsuit})\leq r
$; see the definitions of $K_{\sigma}$ and $R_{\sigma}^{\diamondsuit}$.
Moreover, we have
\[
\bigl\|W^{1/2}K_{\sigma}\bigr\|_2^2=
\Biggl\|W^{1/2}\sum_{i,j=1}^m
\bigl(R_{\sigma
}^{\diamondsuit}\bigr)_{ij} (\phi_i
\otimes\phi_j) \Biggr\|_2^2= \sum
_{i,j=1}^l \lambda_i
\bigl(R_{\sigma}^{\diamondsuit}\bigr)_{ij}^2\leq
\lambda_l \|K_{\sigma}\|_2^2
\]
and
\begin{eqnarray*}
\|K_{\sigma}\|_2^2&=& \Biggl\|\kappa\sum
_{i=1}^{l'}\sum_{j=1}^r
\sigma_{ij}\sum_{k=0}^{[l''/r]-1}
\phi_i\otimes\phi_{l'+rk+j}\\
&&\hspace*{2pt}{} + \kappa\sum
_{i=1}^{r}\sum_{j=1}^{l'}
\sigma_{ji}\sum_{k=0}^{[l''/r]-1}
\phi_{l'+rk+i}\otimes\phi_j \Biggr\|_2^2
\\
&\leq& 2 \kappa^2 l'r \bigl[l''/r
\bigr]\leq\kappa^2 l^2.
\end{eqnarray*}
Therefore,
$
\|W^{1/2}K_{\sigma}\|_{L_2(\Pi^2)}^2 \leq
\lambda_l \kappa^2 \frac{l^2}{m^2},
$
so, we have
%
\begin{equation}
\label{Sobrho} \bigl\|W^{1/2}S_{P_{\sigma}}\bigr\|=\tfrac{1}{16}
\bigl\|W^{1/2}K_{\sigma}\bigr\|_{L_2(\Pi^2)}^2\leq
\rho^2,
\end{equation}
provided that
%
\begin{equation}
\label{kappa3} \kappa\leq\frac{m}{l} \frac{4\rho}{\sqrt{\lambda_l}}.
\end{equation}
We can conclude that, for all $P\in{\cal P}$, $S_{P}\in{\cal
S}_{r,\rho}$
provided that $\kappa$ satisfies conditions~(\ref{kappa1}) and (\ref
{kappa3}).
Since also $|Y|\leq a$, we have that ${\cal P}\subset{\cal P}_{r,\rho,a}$.

Next we check that ${\cal P}$ satisfies the conditions of Proposition
\ref{KBound}.
It is easy to see that, for all $\sigma, \sigma'\in\Lambda' P_{\sigma'}\ll P_{\sigma}$ and
\begin{eqnarray*}
&& K(P_{\sigma}\|P_{\sigma'})\\
&&\qquad={\mathbb E} \biggl( p_{\sigma}
\bigl(X,X'\bigr)\log\frac{p_{\sigma}(X,X')}{p_{\sigma'}(X,X')}+ \bigl(1-p_{\sigma}
\bigl(X,X'\bigr)\bigr)\log\frac{1-p_{\sigma}(X,X')}{1-p_{\sigma'}(X,X')} \biggr).
\end{eqnarray*}
Using the elementary inequality $-\log(1+u)\leq-u+u^2, |u|\leq1/2$
and the fact that $p_{\sigma}(X,X')\in[3/8,5/8],\sigma\in\Lambda$,
we get that
\[
K(P_{\sigma}\|P_{\sigma'})\leq\frac{6}{8^2 a^2}
\|K_{\sigma
}-K_{\sigma'}\|_{L_2(\Pi^2)} \leq\frac{1}{10a^2 m^2}
\|K_{\sigma}-K_{\sigma'}\|_2^2, \sigma,\qquad
\sigma'\in\Lambda'.
\]
A simple computation based on the definition of $K_{\sigma}, K_{\sigma
'}$ easily
yields that
\[
\|K_{\sigma}-K_{\sigma'}\|_2^2\leq8
\kappa^2 l'r \bigl[l''/r
\bigr]\leq8\kappa ^2 l' l''
\leq4 \kappa^2 l^2.
\]
Thus, for the $n$-fold product-measures $P_{\sigma}^{\otimes n},
P_{\sigma'}^{\otimes n}$, we get
\[
K\bigl(P_{\sigma}^{\otimes n}\|P_{\sigma'}^{\otimes n}\bigr)=
nK(P_{\sigma}\|P_{\sigma'})\leq \frac{4n \kappa^2}{10 a^2}
\frac{l^2}{m^2}.
\]
For a fixed $\sigma\in\Lambda'$, this yields
%
\begin{eqnarray}
\label{KLabove} \frac{1}{\operatorname{ card}(\Lambda')-1}\sum_{\sigma'\in\Lambda
'}K
\bigl(P_{\sigma}^{\otimes n}\|P_{\sigma'}^{\otimes n}\bigr)&\leq&
\frac{4 n
\kappa^2}{10 a^2} \frac{l^2}{m^2}
 \leq\frac{1}{10}\frac{l'r}{16}
\nonumber
\\[-8pt]
\\[-8pt]
\nonumber
&\leq&\frac{1}{10} \log\bigl(\operatorname{ card}\bigl(\Lambda'
\bigr)-1\bigr),
\end{eqnarray}
provided that
%
\begin{equation}
\label{kappa2} \kappa\leq\frac{1}{16} a \frac{m}{l}\sqrt{
\frac{rl}{n}}.
\end{equation}

It remains to use (\ref{hemming}) and the definition of kernels
$K_{\sigma}$ to
bound from below the squared distance $\|K_{\sigma}-K_{\sigma'}\|
_{L_2(\Pi^2)}^2$
for $\sigma, \sigma'\in\Lambda', \sigma\neq\sigma'$,
\[
\|K_{\sigma}-K_{\sigma'}\|_{L_2(\Pi^2)}^2
= m^{-2}\|K_{\sigma}-K_{\sigma'}\|_2^2
\geq4 m^{-2} \kappa^2 \frac{l'r}{16}
\bigl[l''/r\bigr] \geq\frac{1}{64}
\kappa^2 \frac{l^2}{m^2}.
\]
Since $S_{P_{\sigma}}=\frac{1}{4}K_{\sigma}$, this implies
that
%
\begin{equation}
\label{L2below} \|S_{P}-S_{P'}\|_{L_2(\Pi^2)}^2
\geq2^{-10}\kappa^2 \frac{l^2}{m^2},\qquad  P,P'
\in{\cal P}, P\neq P'.
\end{equation}

In view of (\ref{kappa1}), (\ref{kappa2}) and (\ref{kappa3}), we now take
\[
\kappa:= \frac{1}{16} a \frac{m}{l}\sqrt{\frac{rl}{n}}
\wedge \frac{m}{l} \frac{4\rho}{\sqrt{\lambda_l}} \wedge
2^{-(1+2/p)}(p-1)^{-1/2} \frac{1}{Q_p(l)}\frac{m}{l}
\frac{a\sqrt {r}}{\sqrt{l}}\frac{1}{m^{2/p}}.
\]
With this choice of $\kappa$,
${\cal P}:=\{P_{\sigma}\dvtx \sigma\in\Lambda'\}\subset{\cal
P}_{r,a,\rho}$.
In view of (\ref{L2below}) and (\ref{KLabove}), we can use
Proposition~\ref{KBound} to get
%
\begin{eqnarray}
\label{knc} &&\inf_{\hat S}\sup_{P\in{\cal P}_{r,a,\rho}}{\mathbb
P}_P \bigl\{ \|\hat S-S_P\|_{L_2(\Pi^2)}^2
\geq c_1\delta_n \bigr\}
\nonumber
\\[-8pt]
\\[-8pt]
\nonumber
&&\qquad\geq \inf_{\hat S}
\sup_{P\in{\cal P}}{\mathbb P}_P \bigl\{ \|\hat
S-S_P\|_{L_2(\Pi^2)}^2\geq c_1
\delta_n \bigr\} \geq c_2,
\end{eqnarray}
where
$
\delta_n:= \frac{a^2rl}{n}
\wedge \frac{\rho^2}{\lambda_l}
\wedge \frac{1}{p-1}\frac{1}{Q_p^2(l)}\frac{a^2 r}{l}\frac{1}{m^{4/p}}
$
and $c_1, c_2>0$ are constants.

In the case when $r>l''$,
bound (\ref{knc}) still holds with
\[
\delta_n:= \frac{a^2 l^2}{n} \wedge \frac{\rho^2}{\lambda_l}
\wedge \frac{1}{p-1}\frac{a^2}{Q_p^2(l)}\frac{1}{m^{4/p}}.
\]
The proof is an easy modification of the argument in the case when
$r\leq l''$.
For $r>l''$, the construction becomes simpler: namely, we define
\[
R^{\flat}:=\pmatrix{ O_{l',l'}
& R & O_{l',m-l}
\vspace*{2pt}\cr
R^{T} & O_{l'',l''} & O_{l'',m-l}
\vspace*{2pt}\cr
O_{m-l,l'} & O_{m-l,l''} & O_{m-l,m-l}},
\]
where $R\in{\cal R}_{l',l''}$, and, based on this,
redefine kernels $K_{\sigma}, \sigma\in\{-1,1\}^{l'\times l''}$.
The proof then goes through with minor simplifications.

Thus, in both cases $r> l''$ and $r\leq l''$, (\ref{knc}) holds with
\[
\delta_n=\delta_n(l):= \frac{a^2(r\wedge l)l}{n} \wedge
\frac{\rho^2}{\lambda_l} \wedge \frac{1}{p-1}
\frac{1}{Q_p^2(l)}\frac{a^2 (r\wedge l)}{l}\frac{1}{m^{4/p}}.
\]
This is true under the assumption that $l\geq32$.
Note also that $Q_p(l)\leq \max_{1\leq j\leq m}\|\bar\phi_j\|
_{L_p(\Pi)}^2=Q_p$.
Thus, we can replace $Q_p^2(l)$ by the upper bound $Q_p^2$ in the definition
of $\delta_n(l)$.

We can now choose $l\in\{32,\ldots, m\}$ that maximizes $\delta_n(l)$
to get bound (\ref{knc}) with $\delta_n:=\min_{32\leq l\leq m}
\delta_n(l)$.
This completes the proof in the case when $k_0\geq32$ and $l_0=32$.
If $k_0<32$, it is easy to use the condition $\lambda_{l+1}\leq
c\lambda_l, l\geq k_0$ and to show that
$
\min_{32\leq l\leq m} \delta_n(l)\leq c'\min_{k_0\leq l\leq m}
\delta_n(l),
$
where $c'$ is a constant depending only on $c$.
This completes the proof in the remaining case.
\end{pf*}

\begin{pf*}{Proof of Theorem~\ref{lowerboundsparse}} The only modification of
the previous proof is to replace bound (\ref{lr-d})\vspace*{2pt} by
$
(\sum_{k=0}^{[l''/r]-1}\phi_{l'+rk+j}(v) )^2\leq
d\sum_{k=0}^{[l''/r]-1}\phi_{l'+rk+j}^2(v).
$
Then, the outcome of the next several lines of the proof is that
${\mathbb P}\{\xi\geq a\}\leq1/4$ provided that [instead of (\ref{kappa1})]
\[
\kappa\leq 2^{-(1+2/p)}(p-1)^{-1/2} \frac{1}{Q_p(l)}
\frac{m}{l}\frac{a}{\sqrt{d}} \frac{1}{m^{2/p}}.
\]
As a result, at the end of the proof, we get that (\ref{knc})
holds with
\[
\delta_n=\delta_n(l):= \frac{a^2(r\wedge l)l}{n} \wedge
\frac{\rho^2}{\lambda_l} \wedge \frac{1}{p-1}
\frac{1}{Q_p^2(l)}\frac{a^2}{d}\frac{1}{m^{4/p}}.
\]
It remains to observe that $Q_p(l)\leq\frac{m}{l}$,
which follows from the fact that
\[
\sum_{j=1}^l\phi_j^2(v)=
\sum_{j=1}^l \langle\phi_j,
e_v\rangle^2 \leq \sum_{j=1}^m
\langle\phi_j, e_v\rangle^2 =1,\qquad v\in V,
\]
and to take $p=\log m$
to complete the proof.
\end{pf*}\vfill\eject

\section{Least squares estimators with nonconvex penalties}\label
{sectLSnonconvex}

In this section, we derive upper bounds on the squared $L_2(\Pi
^2)$-error of the following least squares estimator of the target
matrix $S_{\ast}$:
%
\begin{equation}
\label{leastsquare} \hat S_l:=\hat S_{r,l,a}:=\mathop{\operatorname{argmin}}_{S\in\bar{\cal S}_{r}(l;a)} \frac{1}{n}\sum_{j=1}^n
\bigl(Y_j-S\bigl(X_j,X_j'
\bigr)\bigr)^2,
\end{equation}
where $\bar{\cal S}_{r}(l;a):=\{S^a\dvtx S\in{\cal S}_r(l;a)\}, l=1,\ldots, m$,
\[
{\cal S}_{r}(l;a):= \Biggl\{S\dvtx S\in{\cal
S}_V, \operatorname{ rank}(S)\leq r, \|S\|_{L_2(\Pi^2)}\leq a, S=
\sum_{i,j=1}^l s_{ij}(
\phi_i\otimes\phi_j) \Biggr\}.
\]
Here $S^a$ denotes a truncation
of kernel $S\dvtx S^a(u,v)=S(u,v)$ if $|S(u,v)|\leq a$,
$S^a(u,v)=a$ if $S(u,v)>a$ and $S^a(u,v)=-a$ if $S(u,v)<-a$.
Note that the kernels in the class ${\cal S}_{r}(l;a)$ are symmetric
and $\operatorname{ rank}(S)\leq r\wedge l, S\in{\cal S}_{r}(l;a)$.
Note also that the sets ${\cal S}_{r}(l;a)$, $\bar{\cal S}_{r}(l;a)$
and optimization problem (\ref{leastsquare}) are not convex.
We will prove the following result under the assumption that $|Y|\leq a$
a.s. Recall the definition of the class of kernels ${\cal S}_{r,\rho}$
in Section~\ref{sectlowerbound}.

\begin{theorem}
\label{LSupth}
There exist constants $C>0, A>0$ such that, for all $t>0$, with probability
at least $1-e^{-t}$,
%
\begin{eqnarray}
\label{LSoracul} \|\hat S_l-S_{\ast}\|_{L_2(\Pi^2)}^2
&\leq&2\inf_{S\in\bar{\cal
S}_{r}(l;a)}\|S-S_{\ast}\|_{L_2(\Pi^2)}^2
\nonumber
\\[-8pt]
\\[-8pt]
\nonumber
&&{}+C
\biggl(\frac{a^2 (r\wedge l)l}{n}\log \biggl(\frac{A n m}{(r\wedge
l)l} \biggr)+
\frac{a^2 t}{n} \biggr).
\end{eqnarray}
In particular, for some constants $C,A>0$, for
$S_{\ast}\in{\cal S}_{r,\rho}$ and for all $t>0$,
with probability
at least $1-e^{-t}$,
%
\begin{equation}
\label{LSup} \|\hat S_l-S_{\ast}\|_{L_2(\Pi^2)}^2
\leq C \biggl[\frac{a^2 (r\wedge l)l}{n}\log \biggl(\frac{A n m}{(r\wedge
l)l} \biggr)\vee
\frac{\rho^2}{\lambda_{l+1}}\vee\frac{a^2
t}{n} \biggr].
\end{equation}
\end{theorem}

\begin{pf}
Without loss of generality, assume that $a=1;$ this would imply
the general case by a simple rescaling of the problem. We will use a
version of well-known bounds for least squares estimators over
uniformly bounded function classes
in terms of Rademacher complexities. Specifically, consider the
following least squares estimator:
$
\hat g:= \operatorname{ argmin}_{g\in{\cal G}}n^{-1}\sum_{j=1}^n (Y_j-g(X_j))^2,
$
where $(X_1,Y_1), \ldots, (X_n, Y_n)$ are i.i.d. copies of a random
couple $(X,Y)$
in $T\times{\mathbb R}$, $(T,{\cal T})$ being a measurable space,
$|Y|\leq1$ a.s., ${\cal G}$ being a class of measurable functions on $T$
uniformly bounded by $1$. The goal is to estimate the regression function
$g_{\ast}(x):={\mathbb E}(Y|X=x)$. Define localized Rademacher complexity
$
\psi_n(\delta):=
{\mathbb E}\sup_{g_1,g_2\in{\cal G},\|g_1-g_2\|_{L_2(\Pi)}^2\leq
\delta}|R_n(g_1-g_2)|,
$
where $\Pi$ is the distribution of $X$ and
$
R_n(g):=n^{-1}\sum_{j=1}^n \eps_j g(X_j)
$
is the Rademacher process, $\{\eps_j\}$ being a sequence of i.i.d.
Rademacher random variables
independent of $\{X_j\}$. Denote
$
\psi_n^{\flat}(\delta):=\sup_{\sigma\geq\delta}\frac{\psi
_n(\sigma)}{\sigma}
$
and
$
\psi_n^{\sharp}(\eps):= \inf\{\delta>0\dvtx \psi_n^{\flat}(\delta
)\leq\eps\}.
$
The next result easily follows from Theorem 5.2 in~\cite{Ko27}:

\begin{proposition}
\label{regr}
There exist constants $c_1, c_2>0$ such that, for all $t>0$, with
probability at least $1-e^{-t}$,
\[
\|\hat g-g_{\ast}\|_{L_2(\Pi)}^2\leq 2 \inf
_{g\in{\cal G}}\|g-g_{\ast}\|_{L_2(\Pi)}^2 +
c_1 \biggl(\psi_n^{\sharp}(c_2)+
\frac{t}{n} \biggr).
\]
\end{proposition}

We will apply this proposition to prove Theorem~\ref{LSupth}.
In what follows in the proof, denote $\hat S:=\hat S_l$.
In our case, $T=V\times V$, $(X,X')$ plays the role of $X$, and $\Pi
^2$ plays the role
of $\Pi$. Let ${\cal G}:=\bar{\cal S}_{r}(l;1)$, $g_{\ast}=S_{\ast
}$ and $\hat g=\hat S$.
First, we need to upper bound the Rademacher
complexity $\psi_n(\delta)$ for the class ${\cal G}$.
Let
${\mathbb S}_{r,m}(R)$ be the set of all symmetric $m\times m$ matrices $S$
with $\operatorname{ rank}(S)\leq r$ and $\|S\|_2\leq R$. The $\eps$-covering
number $N({\mathbb S}_{r,m}(R);\|\cdot\|_2;\eps)$ of the set
${\mathbb S}_{r,m}(R)$ with respect to the Hilbert--Schmidt distance
(i.e., the minimal number of balls of radius $\eps$ needed to cover
this set) can be bounded as follows:
%
\begin{equation}
\label{cover2} N\bigl({\mathbb S}_{r,m}(R);\|\cdot\|_2;
\eps\bigr)\leq \biggl(\frac{18
R}{\eps} \biggr)^{(m+1)r}.
\end{equation}
Such bounds are well known
(see, e.g.,~\cite{Ko27}, Lemma 9.3 and references therein; the proof
of this lemma
can be easily modified to obtain (\ref{cover2})).
Bound (\ref{cover2}) will be used to control the covering numbers of
the set of kernels ${\cal S}_{r}(l;1)$. This set can be easily
identified with a subset
of the set ${\mathbb S}_{r\wedge l,l}(m)$ [since kernels $S\in{\cal
S}_{r}(l;1)$
can be viewed as symmetric $l\times l$ matrices of rank at most
$r\wedge l$ with $\|S\|_{L_2(\Pi^2)}\leq1$ and $\|S\|_2=m\|S\|
_{L_2(\Pi^2)}\leq m$]. Therefore, we get the following bound:
$
N({\cal S}_{r}(l;1);\|\cdot\|_2;\eps) \leq (\frac{18 m }{\eps
} )^{(l+1)(r\wedge l)}.
$
Since $\|S_1^1-S_2^1\|_2^2\leq\|S_1-S_2\|_2^2$ (truncation of the entries
reduces the Hilbert--Schmidt distance),
we also have
$
N(\bar{\cal S}_{r}(l;1);\|\cdot\|_2;\eps) \leq (\frac{18 m
}{\eps} )^{(l+1)(r\wedge l)}.
$
Since $\|E_{X_j,X_j'}\|_2\leq1$,
$
\|S_1-S_2\|_{L_2(\Pi_n)}^2 = n^{-1}\sum_{j=1}^n \langle S_1-S_2,
E_{X_j,X_j'}\rangle^2
\leq\|S_1-S_2\|_2^2.
$
Therefore, we get the following bound on the $L_2(\Pi_n)$-covering
numbers of the set $\bar{\cal S}_{r}(l;1)\dvtx
N(\bar{\cal S}_{r}(l;1);L_2(\Pi_n);\eps) \leq (\frac{18 m
}{\eps} )^{(l+1)(r\wedge l)}.
$
Here $\Pi_n$ denotes the empirical distribution based on observations
$(X_1,X_1'),\ldots, (X_n,X_n')$.
The last bound allows us to use inequality (3.17) in~\cite{Ko27} to
control the localized
Rademacher complexity $\psi_n(\delta)$ of the class ${\cal G}$ as follows:
%
\begin{eqnarray}
\label{radem}\quad  \psi_n(\delta)&=&
{\mathbb E}\sup_{S_1,S_2\in\bar{\cal S}_{r}(l;1),\|S_1-S_2\|
_{L_2(\Pi^2)}^2\leq\delta}\Biggl |n^{-1}\sum
_{j=1}^n \eps_j
\bigl(S_1\bigl(X_j,X_j'
\bigr)-S_2\bigl(X_j,X_j'\bigr)
\bigr) \Biggr|
\nonumber
\\[-8pt]
\\[-8pt]
\nonumber
&\leq&
C_1 \biggl[\sqrt{\frac{\delta l(r\wedge l)}{n}}\sqrt{\log \biggl(
\frac{Am}{\sqrt{\delta}} \biggr)} \vee\frac{l(r\wedge l)}{n}\log \biggl(\frac{Am}{\sqrt{\delta
}}
\biggr) \biggr]
\end{eqnarray}
with some constant $A,C_1>0$. This easily yields
$
\psi_n^{\sharp}(c_2)\leq
C_2 \frac{(r\wedge l)l}{n}\log (\frac{Anm}{(r\wedge l)l} )
$
with some constants $A,C_2>0$. Proposition~\ref{regr} now implies
bound (\ref{LSoracul}).

To prove bound (\ref{LSup}), it is enough to observe that, for
$S_{\ast}\in{\cal S}_{r,\rho}$,
%
\begin{equation}
\label{appro} \inf_{S\in\bar{\cal S}_{r}(l;1)}\|S-S_{\ast}
\|_{L_2(\Pi^2)}^2 \leq\frac{2\rho^2}{\lambda_{l+1}}.
\end{equation}
Indeed, since $S_{\ast}\in{\cal S}_{r,\rho}$, we can approximate
this kernel
by
$
S_{l}:=\break\sum_{i,j=1}^l \langle S_{\ast}\phi_i, \phi_j\rangle(\phi
_i\otimes\phi_j).
$
For the error of this approximation, we have
\begin{eqnarray*}
&&\|S_l-S_{\ast}\|_{L_2(\Pi^2)}^2\\[-2pt]
&&\quad= m^{-2}
\|S_l-S_{\ast}\|_2^2=
m^{-2}\sum_{i\vee j>l} \langle S_{\ast}
\phi_i,\phi_j\rangle^2 \\[-2pt]
&&\quad\leq
m^{-2}\frac{1}{\lambda_{l+1}} \sum_{i>l}\sum
_{j=1}^m \lambda _i\langle
S_{\ast}\phi_i,\phi_j\rangle^2
+m^{-2}\frac{1}{\lambda_{l+1}} \sum_{i=1}^m
\sum_{j>l} \lambda _j\langle
S_{\ast}\phi_i,\phi_j\rangle^2\leq
\frac{2\rho
^2}{\lambda_{l+1}},
\end{eqnarray*}
which implies
$
\|S_l^1-S_{\ast}\|_{L_2(\Pi)}^2 \leq\|S_{l}-S_{\ast}\|_{L_2(\Pi^2)}^2
\leq\frac{2\rho^2}{\lambda_{l+1}}
$
(since the entries of matrix $S_{\ast}$ are bounded by $1$ and truncation
of the entries reduces the Hilbert--Schmidt distance).
We also have $\operatorname{ rank}(S_{l})\leq\operatorname{ rank}(S_{\ast})\leq r$ and
\[
\|S_l\|_{L_2(\Pi^2)} = m^{-1} \|S_l
\|_2\leq m^{-1}\|S_{\ast}\|_2=\|
S_{\ast}\|_{L_2(\Pi^2)} \leq\|S_{\ast}\|_{L_{\infty}}\leq1.
\]
Therefore, $S_l^1 \in\bar{\cal S}_r(l;1)$ and bound (\ref{appro}) follows.
Bound (\ref{LSup}) is a consequence of~(\ref{LSoracul}) and (\ref{appro}).
\end{pf}

\begin{remark*}
Note that, in the case when the basis of eigenfunctions
$\{\phi_j\}$
coincides with the canonical basis of space ${\mathbb R}^V$, the following
bound holds trivially:
%
\begin{equation}
\label{LStrivial} \|\hat S_l-S_{\ast}\|_{L_2(\Pi^2)}^2
\leq\frac{4a^2l^2}{m^2}+ \frac
{2\rho^2}{\lambda_{l+1}}.
\end{equation}
This follows from the fact that the entries of both matrices
$\hat S_l$ and $S_l$ are bounded by $a$, and their nonzero entries are only
in the first $l$ rows and the first $l$ columns, so,
$
\|\hat S_l-S_l\|_{L_2(\Pi^2)}^2 \leq\frac{4a^2l^2}{m^2}.
$
Combining this with (\ref{LSup}) and minimizing the resulting bound
with respect
to $l$ yields the following upper bound (up to a constant) that holds
for the optimal choice of $l$:
\[
\min_{1\leq l\leq m} \biggl[ \biggl(\frac{a^2 (r\wedge l)l}{n}\log \biggl(
\frac{A n m}{(r\wedge
l)l} \biggr)\wedge\frac{a^2l^2}{m^2} \biggr)\vee
\frac{\rho
^2}{\lambda_{l+1}} \biggr]\vee\frac{a^2 t}{n}.
\]
It is not hard to check that, typically, this expression is of the same order
(up to log factors) as the lower bound of Theorem~\ref
{lowerboundsparse} for
$d=1$.
\end{remark*}

Next we consider a penalized version of least squares estimator which
is adaptive
to unknown parameters of the problem (such as the rank of the target matrix\vadjust{\goodbreak}
and the optimal value of parameter $l$ which minimizes the error bound
of Theorem~\ref{LSupth}). We still assume that $|Y|\leq a$
a.s. for some known constant $a>0$.
Define
%
\begin{eqnarray}
\label{penLS}  (\hat r, \hat l)&:=& \mathop{\operatorname{ argmin}}_{r,l=1,\ldots, m}
\Biggl\{n^{-1}\sum_{j=1}^n
\bigl(Y_j-\hat S_{r,l,a}\bigl(X_j,X_j'
\bigr)\bigr)^2
\nonumber
\\[-8pt]
\\[-8pt]
\nonumber
&&\hspace*{45pt}{}+K\frac{a^2 (r\wedge l)l}{n}\log \biggl(\frac{A n m}{(r\wedge
l)l} \biggr) \Biggr\}
\end{eqnarray}
and let $\hat S:=\hat S_{\hat r, \hat l, a}$.
Here $K>0$ and $A>0$ are fixed constants.

The following theorem provides an oracle inequality for the estimator
$\hat S$.

\begin{theorem}
\label{regradapt}
There exists a choice of constants $K>0$, $A>0$ in (\ref{penLS}) and
$C>0$ in the inequality below such that for all $t>0$ with probability
at least $1-e^{-t}$
%
\begin{eqnarray}
\label{LSoracul01} &&
\|\hat S-S_{\ast}\|_{L_2(\Pi^2)}^2
\nonumber
\\
&&\qquad\leq2\min_{1\leq r\leq m,1\leq
l\leq m} \biggl[\inf_{S\in\bar{\cal S}_{r}(l;a)}
\|S-S_{\ast}\| _{L_2(\Pi^2)}^2\\
&&\hspace*{100pt}{}+
C \biggl(\frac{a^2 (r\wedge l)l}{n}\log \biggl(\frac{A n m}{(r\wedge
l)l} \biggr)+
\frac{a^2 (t+\log m)}{n} \biggr) \biggr].\nonumber
\end{eqnarray}
\end{theorem}

\begin{pf} As in the proof of the previous theorem, we can assume that
$a=1;$ the general case follows by rescaling.
We will use oracle inequalities in abstract penalized empirical risk
minimization problems; see~\cite{Ko27}, Theorem~6.5.
We only sketch the proof here skipping the details that are standard.
As in the proof of Theorem~\ref{LSupth}, first consider
i.i.d. copies $(X_1,Y_1), \ldots, (X_n, Y_n)$ of a random couple
$(X,Y)$ in $T\times{\mathbb R}$, where $(T,{\cal T})$ is a measurable
space and $|Y|\leq1$ a.s.
Let $\{{\cal G}_k\dvtx k\in I\}$ be a finite family of classes of measurable
functions from $T$ into $[-1,1]$. Consider the corresponding family of
least squares estimators
$
\hat g_k:= \operatorname{ argmin}_{g\in{\cal G}_k}n^{-1}\sum_{j=1}^n
(Y_j-g(X_j))^2, k\in I.
$
Suppose the following upper bounds on localized Rademacher complexities for
classes ${\cal G}_k, k\in I$ hold:
$
{\mathbb E}\sup_{g_1,g_2\in{\cal G}_k,\|g_1-g_2\|_{L_2(\Pi)}^2\leq
\delta}|R_n(g_1-g_2)|
\leq\psi_{n,k}(\delta), \delta>0,
$
where $\psi_{n,k}$ are nondecreasing functions of $\delta$ that
do not depend on the distribution of $(X,Y)$. Let
%
\begin{equation}
\label{regpe} \hat k:= \mathop{\operatorname{argmin}}_{k\in I}
\Biggl[n^{-1}\sum_{j=1}^n
\bigl(Y_j-\hat g_k(X_j)\bigr)^2
+ K \biggl(\psi_{n,k}^{\sharp}(c_1)+
\frac{t_k}{n} \biggr) \Biggr],
\end{equation}
and $K,c_1$ are constants and $\{t_k, k\in I\}$ are positive
numbers. Define the following penalized least squares estimator of the
regression function \mbox{$g_{\ast}\dvtx \hat g:=\hat g_{\hat k}$}.

The next result is well known; it can be deduced, for instance, from
Theorem~6.5 in~\cite{Ko27}.

\begin{proposition}
\label{regror}
There exists constants $K,c_1>0$ in the definition (\ref{regpe}) of
$\hat k$
and a constant $K_1>0$
such that, for all $t_k>0$, with
probability at least $1-\sum_{k\in I}e^{-{t_k}}$
\[
\|\hat g-g_{\ast}\|_{L_2(\Pi)}^2\leq 2 \inf
_{k\in I} \biggl[\inf_{g\in{\cal G}_k}\|g-g_{\ast}
\|_{L_2(\Pi)}^2 + K_1 \biggl(\psi_{n,k}^{\sharp}(c)+
\frac{t_k}{n} \biggr) \biggr].
\]
\end{proposition}

We apply this result to the estimator $\hat S=\hat S_{\hat r, \hat l,1}$,
where $(\hat r, \hat l)$ is defined by (\ref{penLS}) (with $a=1$).
In this case, $T=V\times V$, $(X,X')$ plays the role of $X$, $g_{\ast
}=S_{\ast}$,
$I=\{(r,l)\dvtx 1\leq r,l\leq m\}$, ${\cal G}_{r,l}=\bar{\cal S}_{r}(l;1)$.
In view of (\ref{radem}), we can use the following bounds on localized
Rademacher complexities for these function classes:
\[
\psi_{n,r,l}(\delta):=C_1 \biggl[\sqrt{\frac{\delta l(r\wedge
l)}{n}}
\sqrt{\log \biggl(\frac{Am}{\sqrt{\delta}} \biggr)} \vee\frac{l(r\wedge l)}{n}\log \biggl(
\frac{Am}{\sqrt{\delta
}} \biggr) \biggr]
\]
with some constant $C_1$, and we have
$
\psi_{n,r,l}^{\sharp}(c_1)\leq
C_2 \frac{(r\wedge l)l}{n}\log (\frac{Anm}{(r\wedge l)l} )
$
with some constant $C_2>0$. Define $t_{r,l}:=t+2\log m, (r,l)\in I$.
This yields the bound
$
\sum_{(r,l)\in I}e^{-t_{r,l}}\leq e^{-t}.
$
These considerations and Proposition~\ref{regror} imply the claim
of the theorem.
\end{pf}

It follows from Theorem~\ref{regradapt} that, for some constant $C>0$ and
for all $t>0$,
%
\begin{equation}\label{LSupad} \sup_{P\in{\cal P}_{r,\rho,a}}{\mathbb P}_P \biggl
\{\|\hat S-S_P\| _{L_2(\Pi^2)}^2\geq C \biggl(
\Delta_n(r,\rho,a)\vee\frac{a^2 t}{n} \biggr) \biggr\} \leq
e^{-t},
\end{equation}
where
$
\Delta_n(r,\rho,a):=
\min_{1\leq l\leq m} [\frac{a^2 (r\wedge l)l}{n}\log
(\frac{An m}{(r\wedge l)l} )\vee\frac{\rho^2}{\lambda
_{l+1}} ].
$
Denoting
\[
\tilde l:= \min \biggl\{l=1,\ldots, m\dvtx (r\vee l)l \lambda_{l+1}
\log \biggl(\frac{An m}{(r\wedge l)l} \biggr)\geq\frac{\rho^2
n}{a^2} \biggr\},
\]
it is easy to see that
$
\Delta_n(r,\rho,a)=
\frac{a^2(r\wedge\tilde l)\tilde l}{n}\log (\frac{An
m}{(r\wedge\tilde l)\tilde l} ) \vee\frac{\rho^2}{\lambda
_{\tilde l}}.
$

\begin{example*} Suppose that, for some $\beta>1/2$, $\lambda_l\asymp
l^{2\beta},
l=1, \ldots, m$. Under this assumption, it is easy to show that the
upper bound
on the squared $L_2(\Pi^2)$-error of the estimator $\hat S$ is of the order
\begin{eqnarray*}
&&\biggl( \biggl(\frac{a^2 \rho^{1/\beta}r}{n}\log\frac{Anm}{r}
\biggr)^{2\beta/(2\beta+1)} \wedge \biggl(\frac{a^2\rho^{2/\beta}\log(Anm)}{n}
\biggr)^{\beta/\beta+1} \\
&&{}\hspace*{170pt}\qquad\wedge\frac{a^2 rm \log(Anm)}{n} \biggr)\vee
\frac{a^2 t}{n}
\end{eqnarray*}
(in fact, the log factors can be written in a slightly better, but
more complicated way). Up to the log factors, this is the same error rate
as in the lower bounds of Section~\ref{sectlowerbound}; see (\ref
{minmaxlower}).
\end{example*}

\section{Least squares with convex penalization:
Combining nuclear norm and squared Sobolev norm}\label{sectLSconvex}

Our main goal in this section is to study the following penalized least
squares estimator with a combination of two convex penalties:
%
\begin{equation}
\label{LSnucSob} \hat S_{\eps, \bar\eps}:= \mathop{\operatorname{ argmin}}_{S\in{\mathbb D}}
\Biggl[ \frac{1}{n}\sum_{j=1}^n
\bigl(Y_j -S\bigl(X_j,X_j'
\bigr)\bigr)^2+ \eps\|S\|_1 + \bar\eps\bigl\|W^{1/2}S
\bigr\|_{L_2(\Pi^2)}^2 \Biggr],\hspace*{-35pt}
\end{equation}
where ${\mathbb D}\subset{\cal S}_V$ is a closed convex set of
symmetric kernels such that, for all $S\in{\mathbb D}$,
$
\|S\|_{L_{\infty}}:=\max_{u,v\in V}|S(u,v)|\leq a,
$
and $\eps, \bar{\eps}>0$ are regularization parameters.
The first penalty involved in (\ref{LSnucSob}) is based on the
nuclear norm $\|S\|_1$, and it is used to ``promote'' low-rank
solutions. The second penalty is based on a ``Sobolev type norm'' $\|
W^{1/2}S\|_{L_2(\Pi^2)}^2$. It is used to ``promote'' the smoothness
of the solution on the graph.

We will derive an upper bound on the error
$\|\hat S_{\eps,\bar\eps}-S_{\ast}\|_{L_2(\Pi^2)}^2$ of estimator
$\hat S_{\eps,\bar\eps}$ in terms of spectral characteristics of the
target kernel $S_{\ast}$ and matrix~$W$.
As before, $W$ is a nonnegatively definite symmetric kernel
with spectral representation
$
W=\sum_{k=1}^m \lambda_k (\phi_k \otimes\phi_k),
$
where $0\leq\lambda_1 \leq\cdots\leq\lambda_m$ are the eigenvalues
of $W$ repeated with their multiplicities and $\phi_1,\ldots, \phi_m$
are the corresponding orthonormal eigenfunctions.
We will also use the decomposition of identity
associated with $W$:
$
E(\lambda):=
\sum_{\lambda_j\leq\lambda}(\phi_j\otimes\phi_j), \lambda\geq0.
$
Clearly, $\lambda\mapsto E(\lambda)$ is a nondecreasing
projector-valued function. Despite the fact that the eigenfunctions $\{
\phi_k\}$ are not
uniquely defined in the case when $W$ has multiple eigenvalues, the
decomposition
of identity $\{E(\lambda), \lambda\geq0\}$ is uniquely defined (in
fact, it can be
rewritten in terms of spectral projectors of $W$).
The distribution of the eigenvalues of $W$ is characterized by
the following \emph{spectral function}:
\[
F(\lambda):=\operatorname{ tr}\bigl(E(\lambda)\bigr)=\bigl\|E(\lambda)
\bigr\|_2^2= \sum_{j=1}^m
I(\lambda_j \leq\lambda), \qquad\lambda\geq0.
\]
Denote $k_0: = F(0)+1$ (in other words, $k_0$ is the smallest $k$ such
that $\lambda_k> 0$).
It was assumed in the \hyperref[intro]{Introduction} that there exists a constant $c\geq
1$ such that
$\lambda_{k+1}\leq c \lambda_k$ for all $k\geq k_0$.

In what follows, we use a regularized majorant of spectral function $F$.
Let $\bar F\dvtx {\mathbb R}_+\mapsto{\mathbb R}_+$ be a
nondecreasing function such that
$F(\lambda)\leq\bar F(\lambda), \lambda\geq0$, the function
$\lambda\mapsto\frac{\bar F(\lambda)}{\lambda}$ is nonincreasing
and, for some $\gamma\in(0,1)$,
\[
\int_{\lambda}^{\infty}\frac{\bar F(s)}{s^2}\,ds\leq
\frac{1}{\gamma}\frac{\bar F(\lambda)}{\lambda}, \qquad\lambda>0.
\]
Without loss of
generality, we assume in what follows that $\bar F(\lambda)=m, \lambda
\geq\lambda_m$ [otherwise, one can take the function $\bar F(\lambda
)\wedge m$ instead].
The conditions on $\bar F$ are satisfied if for some\vadjust{\goodbreak} $\gamma\in(0,1)$,
the function
$\frac{\bar F(\lambda)}{\lambda^{1-\gamma}}$ is nonincreasing:
in this case, $\frac{\bar F(\lambda)}{\lambda}$ is also
nonincreasing and
\[
\int_{\lambda}^{\infty}\frac{\bar F(s)}{s^2}\,ds= \int
_{\lambda}^{\infty}\frac{\bar F(s)}{s^{1-\gamma}}\frac
{ds}{s^{1+\gamma}}
\leq \frac{\bar F(\lambda)}{\lambda^{1-\gamma}} \int_{\lambda}^{\infty}
\frac{ds}{s^{1+\gamma}} =\frac{1}{\gamma}\frac{\bar F(\lambda)}{\lambda}.
\]

Consider a kernel $S\in{\cal S}_V$ (an oracle) with spectral
representation: $S=\sum_{k=1}^r \mu_k (\psi_k \otimes\psi_k)$, where
$r=\operatorname{ rank}(S)\geq1$, $\mu_k$ are nonzero eigenvalues of $S$
(possibly repeated) and $\psi_k$ are the corresponding orthonormal
eigenfunctions. Denote $L=\operatorname{ supp}(S)=\mbox{l.s.}(\psi_1,\ldots, \psi
_r)$. The following \emph{coherence function} will be used
to characterize the relationship between the kernels $S$ and~$W$:
%
\begin{equation}
\label{cohere} \varphi(S;\lambda):= \bigl\langle P_L, E(\lambda)
\bigr\rangle:= \sum_{\lambda_j\leq\lambda}\|P_L
\phi_j\|^2,\qquad \lambda\geq0.
\end{equation}
It is immediate from this definition that
$
\varphi(S,\lambda)\leq F(\lambda)\leq\bar F(\lambda),
\lambda\geq0.
$
Note also that $\varphi(S;\lambda)$ is a nondecreasing function
of $\lambda$ and
$
\varphi(S,\lambda)=\break\sum_{j=1}^m\|P_L\phi_j\|^2=r, \lambda\geq
\lambda_m
$
[for $\lambda<\lambda_m$, $\varphi(S;\lambda)$ can be interpreted
as a ``partial rank'' of $S$]. As in the case of spectral function $S$,
we need
a regularized majorant for the coherence function $\varphi(S;\lambda)$.
Denote by $\Psi=\Psi_{S,W}$ the set of all nondecreasing functions
$\varphi\dvtx {\mathbb R}_+\mapsto{\mathbb R}_+$ such that $\lambda
\mapsto\frac{\varphi(\lambda)}{\bar F(\lambda)}$ is nonincreasing
and $\varphi(S;\lambda)\leq\varphi(\lambda), \lambda\geq0$. It
is easy to see that the class of functions $\Psi_{S,W}$ contains
the smallest function (uniformly in $\lambda\geq0$) that will be denoted
by $\bar\varphi(S;\lambda)$ and it is given by the following expression:
\[
\bar\varphi(S;\lambda):= \sup_{\sigma\leq\lambda} \bar F(\sigma )\sup
_{\sigma'\geq\sigma} \frac{\varphi(S;\sigma')}{\bar F(\sigma')}.
\]
It easily follows from this definition that $\bar\varphi(S,\lambda
)=r, \lambda\geq\lambda_m$. Note that since the function $\frac
{\bar\varphi(S,\lambda)}{\bar F(\lambda)}$ is nonincreasing and it
is equal to $\frac{r}{m}$ for $\lambda\geq\lambda_m$, we have
%
\begin{equation}
\label{lowvarphi} \bar\varphi(S;\lambda)\geq\frac{r}{m}\bar F(\lambda)
\geq\frac
{r}{m}F(\lambda), \qquad\lambda\geq0.
\end{equation}

Given $t>0$, $\tilde\lambda\in(0,\lambda_{k_0}]$, let
$t_{n,m}:=
t+
3\log (2\log_2 n+\frac{1}{2}\log_2\frac{\lambda_m}{\tilde
\lambda}+2 ).
$
Suppose that, for some $D>0$,
%
\begin{equation}
\label{condeps} \eps\geq D a \biggl(\sqrt{\frac{\log(2m)}{nm}}\vee
\frac{\log
(2m)}{n} \biggr).
\end{equation}

\begin{theorem}
\label{main}
There exists constants $C, D$ depending only on $c,\gamma$ such that,
for all $\bar\eps\in[0,\tilde\lambda^{-1}]$
with probability at least $1-e^{-t}$,
%
\begin{eqnarray}
\label{bound_th_main} && \|\hat S_{\eps,\bar\eps}-S_{\ast}
\|_{L_2(\Pi^2)}^2
\nonumber
\\
&&\qquad\leq \inf_{S\in{\mathbb D}} \bigl[
\|S-S_{\ast}\|_{L_2(\Pi^2)}^2+ C m^2
\eps^2 \bar\varphi\bigl(S;\bar\eps^{-1}\bigr)
+\bar\eps\bigl\|W^{1/2}S\bigr\|_{L_2(\Pi^2)}^2 \bigr]
\\
&&\qquad\quad{}+ C\frac{a^2 t_{n,m}}{n}.\nonumber
\end{eqnarray}
\end{theorem}

\begin{remarks*}
(1) Under the additional assumption that $m \log(2m)\leq n$, one can take
$
\eps= Da\sqrt{\frac{\log(2m)}{nm}}$.\vspace*{1pt}
In this case, the main part of the random error term in the right-hand
side of bound (\ref{bound_th_main}) becomes
\begin{eqnarray*}
&& C m^2 \eps^2 \bar\varphi\bigl(S;\bar\eps^{-1}
\bigr) +\bar\eps\bigl\|W^{1/2}S\bigr\|_{L_2(\Pi^2)}^2\\
&&\qquad =
C'\frac{a^2 \bar\varphi(S;\bar\eps^{-1})m\log(2m)}{n}+\bar\eps\bigl\| W^{1/2}S
\bigr\|_{L_2(\Pi^2)}^2.
\end{eqnarray*}

(2) Note also that Theorem~\ref{main} holds in the case when $\bar
\eps=0$.
In this case, our method coincides with nuclear norm penalized least
squares (matrix LASSO) and
$\bar\varphi(S; \bar\eps^{-1})=\operatorname{ rank}(S)$,
so the bound of Theorem~\ref{main} becomes
%
\begin{equation}
\label{worsebound} \|\hat S_{\eps,0}-S_{\ast}
\|_{L_2(\Pi^2)}^2\leq\inf_{S\in{\mathbb
D}} \bigl[
\|S-S_{\ast}\|_{L_2(\Pi^2)}^2+ C m^2
\eps^2 \operatorname{ rank}(S) \bigr]+ C\frac{a^2 t_{n,m}}{n}.\hspace*{-35pt}
\end{equation}
Similar oracle inequalities were proved in~\cite{Ko26}
for a linearized
least squares method with nuclear norm penalty.
\end{remarks*}

Using simple aggregation techniques, it is easy to construct an adaptive
estimator for which the oracle inequality of Theorem~\ref{main} holds
with the optimal value of $\bar\eps$ that minimizes the right-hand
side of the bound. To this end, divide the sample $(X_1,X_1',Y_1),\ldots,
(X_n,X_n',Y_n)$ into two parts,
\begin{eqnarray*}
&&\bigl(X_j,X_j',Y_j\bigr),\qquad
j=1,\ldots, n'\quad \mbox{and}\\
 &&\bigl(X_{n'+j},
X_{n'+j}', Y_{n'+j}\bigr),\qquad j=1,\ldots,
n-n',
\end{eqnarray*}
where $n':=[n/2]+1$. The first part of the sample will be used
to compute the estimators $\hat S_l:=\hat S_{\eps, \bar\eps_l}$,
$\eps_l:=\lambda_{l}^{-1}$, $l=k_0,\ldots, m+1$ [they are defined
by (\ref{LSnucSob}), but they are based only on the first
$n'$ observations]. The second part of the sample is used
for model selection
\[
\hat l:= \mathop{\operatorname{ argmin}}_{l=k_0,\ldots, m+1}\frac{1}{n-n'}\sum
_{j=1}^{n-n'} \bigl(Y_{n'+j}-\hat
S_{l}\bigl(X_{n'+j},X_{n'+j}'\bigr)
\bigr)^2.
\]
Finally, let $\hat S:=\hat S_{\hat l}$.

\begin{theorem}
\label{mainadapt}
Under the assumptions and notation of Theorem~\ref{main},
with probability at least $1-e^{-t}$,
%
\begin{eqnarray}
\label{boundthmainadapt}
\|\hat S-S_{\ast}\|_{L_2(\Pi^2)}^2
&\leq& \inf_{S\in{\mathbb D}} \Bigl[2\|S-S_{\ast}
\|_{L_2(\Pi^2)}^2
\nonumber\\
&&\hspace*{30pt}{}+ C\inf_{\bar\eps\in[0,\lambda_{k_0}^{-1}]}
\bigl(m^2 \eps^2 \bar\varphi\bigl(S;\bar
\eps^{-1}\bigr)
+\bar\eps\bigl\|W^{1/2}S\bigr\|_{L_2(\Pi^2)}^2 \bigr)
\Bigr] \\
&&{}+ C\frac{a^2 (\log(m+1)+t_{n,m})}{n}.\nonumber
\end{eqnarray}
\end{theorem}

\begin{pf} The idea of aggregation result behind this theorem is rather
well known; see~\cite{Massart}, Chapter 8. The proof can be deduced,
for instance, from Proposition~\ref{regr}
used in Section~\ref{sectLSnonconvex}. Specifically, this
proposition has to be applied in the case when ${\cal G}$ is a finite
class of functions bounded by $1$. Let $N:=\operatorname{ card}({\cal G})$. Then,
for some numerical constant $C_1>0$
\[
\psi_n(\delta)\leq C_1 \biggl[\delta\sqrt{
\frac{\log N}{n}}\vee \frac{\log N}{n} \biggr]
\]
(see, e.g.,~\cite{Ko27}, Theorem 3.5), and Proposition~\ref{regr} easily
implies that, for all $t>0$, with probability at least $1-e^{-t}$
%
\begin{equation}
\label{finiteaggregate} \|\hat g-g_{\ast}\|_{L_2(\Pi)}^2
\leq 2 \inf_{g\in{\cal G}}\|g-g_{\ast}\|_{L_2(\Pi)}^2
+ C_2 \frac{\log N+t}{n},
\end{equation}
where $C_2>0$ is a constant. We will assume that $a=1$ (in the general case,
the result would follow by rescaling) and
use bound (\ref{finiteaggregate}), conditionally on the
first part of the sample, in the case when ${\cal G}:=\{\hat g_l\dvtx l=k_0,\ldots, m+1\}$.
Then, given $(X_j,X_j',Y_j), j=1,\ldots, n'$,
with probability
at least $1-e^{-t}$,
%
\begin{equation}
\label{qua-qua} \|\hat S-S_{\ast}\|_{L_2(\Pi^2)}^2\leq 2
\min_{k_0\leq l\leq m+1}\|\hat S_l-S_{\ast}
\|_{L_2(\Pi)}^2 + C_2 \frac{\log(m+1)+t}{n}.
\end{equation}
By Theorem~\ref{main} [with $t$ replaced by $t+\log(m+1)$] and the
union bound, we get that, with probability
at least $1-e^{-t}$, for all $l=k_0,\ldots, m+1$,
%
\begin{eqnarray}
\label{boundthmain} && \|\hat S_l-S_{\ast}
\|_{L_2(\Pi^2)}^2
\nonumber
\\
&&\qquad\leq \inf_{S\in{\mathbb D}} \bigl[
\|S-S_{\ast}\|_{L_2(\Pi^2)}^2+ C_3
m^2 \eps^2\bar\varphi\bigl(S;\bar\eps_l^{-1}
\bigr)
+\bar\eps_l \bigl\|W^{1/2}S\bigr\|_{L_2(\Pi^2)}^2
\bigr]\\
&&\qquad\quad{} + C_3\frac{\log(m+1)+t_{n,m}}{n}\nonumber
\end{eqnarray}
with some constant $C_3>0$. Therefore,
the minimal error of estimators $\hat S_l$, $\min_{k_0\leq l\leq m+1}\|
\hat S_l-S_{\ast}\|_{L_2(\Pi)}^2$, can be bounded with the same probability
by the minimum over $l=k_0, \ldots, m+1$ of the expression in
the right-hand side of (\ref{boundthmain}).
Moreover, using monotonicity of the function $\lambda\mapsto\varphi
(S;\lambda)$
and the condition that $\lambda_{l+1}\leq c\lambda_l, l=k_0,\ldots,
m-1$, it is easy to replace
the minimum over $l$ by the infimum over $\bar\eps$.
Combining the resulting bound with (\ref{qua-qua}) and adjusting the
constants yields the claim.
\end{pf}

Using more sophisticated aggregation methods (e.g., such as the methods
studied in~\cite{GaifasLecue}) it is possible to construct an estimator
$\hat S$ for which the oracle inequality similar to (\ref{boundthmainadapt})
holds with constant $1$ in front of the approximation error term
$\|S-S_{\ast}\|_{L_2(\Pi^2)}^2$.

To understand better the meaning of function $\bar\varphi$ involved
in the statements of Theorems~\ref{main} and~\ref{mainadapt}, it
makes sense to
relate it to the low coherence assumptions discussed in the
\hyperref[intro]{Introduction}.\vadjust{\goodbreak} Indeed, suppose that, for some $\nu=\nu(S)\geq1$,
%
\begin{equation}
\label{lowcoherence} \|P_L \phi_k\|^2 \leq
\frac{\nu r}{m},\qquad k=1,\ldots, m.
\end{equation}
This is a part of standard low coherence assumptions on matrix $S$ with
respect to the orthonormal basis $\{\phi_k\}$; see (\ref{coherA}).
Clearly, it implies that\footnote{Compare (\ref{up_varphi}) with
(\ref{lowvarphi}).}
%
\begin{equation}
\label{up_varphi} \bar\varphi(S;\lambda)\leq\frac{\nu r \bar F(\lambda)}{m}, \qquad\lambda\geq0.
\end{equation}
Suppose
that $n\geq m\log(2m)$ and
$
\eps= Da\sqrt{\frac{\log(2m)}{nm}}.
$
If condition (\ref{up_varphi}) holds for the target kernel $S_{\ast}$
with $r=\operatorname{ rank}(S_{\ast})$ and some $\nu\geq1$, then Theorem~\ref{main}
implies that with probability at least $1-e^{-t}$,
\begin{eqnarray*}
\|\hat S_{\eps,\bar\eps}-S_{\ast}\|_{L_2(\Pi^2)}^2&\leq& C
\frac{a^2\nu r \bar F(\bar\eps^{-1}) \log(2m)}{n} +\bar\eps\bigl\|W^{1/2}S_{\ast}
\bigr\|_{L_2(\Pi^2)}^2\\
&&{}+ C \frac{a^2 t_{n,m}}{n},
\end{eqnarray*}
and Theorem~\ref{mainadapt} implies that with the same probability,
\begin{eqnarray*}
\|\hat S-S_{\ast}\|_{L_2(\Pi^2)}^2&\leq& C\inf
_{\bar\eps\in[0,\lambda_{k_0}^{-1}]} \biggl(\frac{a^2\nu r
\bar F(\bar\eps^{-1}) \log(2m)}{n} +\bar\eps
\bigl\|W^{1/2}S_{\ast}\bigr\|_{L_2(\Pi^2)}^2 \biggr)\\
&&{}+ C
\frac{a^2(\log(m+1)+ t_{n,m})}{n}.
\end{eqnarray*}

\begin{example*} If $\lambda_k \asymp k^{2\beta}$ for some $\beta
>1/2$, then
it is easy to check that $\bar F(\lambda)\asymp\lambda^{1/2\beta}$.
Under the assumption that $\|W^{1/2}S_{\ast}\|_{L_2(\Pi^2)}^2\leq
\rho^2$,
we get the bound
%
\begin{eqnarray}
\label{betterbound} && \|\hat S-S_{\ast}\|_{L_2(\Pi^2)}^2
\nonumber
\\
&&\qquad\leq C \biggl( \biggl( \biggl(\frac{a^2 \rho^{1/\beta}\nu r \log
(2m)}{n} \biggr)^{2\beta/(2\beta+1)} \wedge
\frac{a^2 r m}{n} \biggr)
\\
&&\hspace*{92pt}\qquad\quad{}\vee\frac{a^2 (\log(m+1)+t_{n,m})}{n} \biggr).\nonumber
\end{eqnarray}
Under the following slightly modified version of low coherence
assumption~(\ref{up_varphi}),
%
\begin{equation}
\label{upvarphi} \bar\varphi(S;\lambda)\leq\frac{\nu(r\wedge\bar F(\lambda))\bar
F(\lambda)}{m},\qquad \lambda\geq0,
\end{equation}
one can almost recover upper bounds of Section~\ref{sectLSnonconvex},
\begin{eqnarray*}
\label{betterbound} &&
\|\hat S-S_{\ast}\|_{L_2(\Pi^2)}^2
\\
&&\qquad
\leq C \biggl( \biggl( \biggl(\frac{\nu a^2 \rho^{1/\beta}r \log
(2m)}{n} \biggr)^{2\beta/(2\beta+1)} \wedge
\biggl(\frac{\nu a^2 \rho^{2/\beta} \log(2m)}{n} \biggr)^{\beta
/(\beta+1)}
\\
&&\hspace*{288pt}{}\wedge\frac{a^2 r m}{n} \biggr)\\
&&\hspace*{222pt}{} \vee
\frac{a^2 (\log(m+1)+t_{n,m})}{n} \biggr).
\end{eqnarray*}
The main difference with what was proved in Section \ref
{sectLSnonconvex} is
that now the low coherence constant $\nu$ is involved in the bounds, so
the methods discussed in this section yield correct (up to log factors)
error rates
provided that the target kernel $S_{\ast}$ has ``low coherence'' with respect
to the basis of eigenfunctions of $W$.
\end{example*}

\begin{pf*}{Proof of Theorem~\ref{main}} Bound (\ref{bound_th_main}) will
be proved
for a fixed oracle $S\in{\mathbb D}$ and an arbitrary function
$\varphi\in\Psi_{S,W}$ with
$\varphi(\lambda)=r, \lambda\geq\lambda_m$ instead of $\bar
\varphi$. It then can be applied to the function $\bar\varphi$
(which is the smallest function in $\Psi_{S,W}$).
Without loss of generality, we assume that $a=1$;
the general case then follows by a simple rescaling.
Finally, we will denote $\hat S:=\hat S_{\eps,\bar\eps}$ throughout
the proof.

Define the following orthogonal projectors ${\mathcal P}_L, {\mathcal
P}_L^{\perp}$
in the space ${\cal S}_V$ with Hilbert--Schmidt inner product:
$
\mathcal{P}_L(A):=A-P_{L^\perp}A P_{L^\perp},  \mathcal{P}_L^\perp
(A)=P_{L^\perp}A P_{L^\perp}, \break A\in{\mathcal S}_V.
$
We will use a well known representation of subdifferential of convex
function $S \mapsto\|S\|_1$:
\[
\partial\|S\|_1= \bigl\{ \operatorname{ sign}(S) +
\mathcal{P}_L^\perp(M) \dvtx M \in\mathcal{S}_V,
\| M \| \leq1 \bigr\},
\]
where $L=\operatorname{ supp}(S)$; see~\cite{Ko27}, Appendix A.4 and references
therein. Denote
\[
L_n(S):= \frac{1}{n}\sum_{j=1}^n
\bigl(Y_j -S\bigl(X_j,X_j'
\bigr)\bigr)^2+ \eps\|S\|_1 + \bar\eps\bigl\|W^{1/2}S
\bigr\|_{L_2(\Pi^2)}^2,
\]
so that $\hat S:=\operatorname{ argmin}_{S\in{\mathbb D}}L_n(S)$.
An arbitrary matrix $A \in\partial L_n(\hat{S})$ can be represented as
%
\begin{equation}
\label{Aravno} \qquad A=\frac{2}{n}\sum_{i=1}^n
\hat S\bigl(X_i,X_i'\bigr)
E_{X_i,X_i'} - \frac
{2}{n}\sum_{i=1}^{n}
Y_i E_{X_i,X_i'}+\eps\hat{V}+2\frac{\bar\eps
}{m^2} W \hat{S},
\end{equation}
where $\hat{V} \in\partial\| \hat{S} \|_1$.
Since $\hat{S}$ is a minimizer of $L_n(S)$, there exists a matrix $A
\in\partial L_n(\hat{S})$ such that $-A$ belongs to the normal cone
of $\mathbb{D}$ at the point $\hat{S}$; see~\cite{Aubin},\vadjust{\goodbreak} Chapter 2,
Corollary 6.
This implies that $\langle A, \hat{S} - S \rangle\leq0$ and, in view of
(\ref{Aravno}),
%
\begin{eqnarray}
\label{Aravno1}&& 2P_n\bigl(\hat S (\hat S-S)\bigr)- \Biggl\langle
\frac{2}{n} \sum_{i=1}^{n}
Y_i E_{X_i,X_i'},\hat{S}-S \Biggr\rangle+
\eps\langle\hat{V},\hat{S}-S \rangle
\nonumber
\\[-8pt]
\\[-8pt]
\nonumber
&&\qquad{}+ 2\frac{\bar\eps}{m^2} \langle
W\hat{S},\hat{S}-S \rangle\leq0.
\end{eqnarray}
Here and in what follows $P_n$ denotes the empirical
distribution based on the sample $(X_1,X_1',Y_1),\ldots, (X_n,X_n',Y_n)$.
The corresponding true distribution of $(X,X',Y)$ will be denoted by $P$.
It easily follows from (\ref{Aravno1}) that
\begin{eqnarray*}
&&
2\langle\hat S-S_{\ast}, \hat S-S\rangle_{L_2(P_n)} -2
\langle\Xi,\hat{S}-S \rangle
\\
&&\qquad{}+
\eps\langle\hat{V},\hat{S}-S \rangle+ 2\bar\eps \bigl\langle
W^{1/2}\hat{S},W^{1/2}(\hat{S}-S) \bigr\rangle_{L_2(\Pi^2)}
\leq0,
\end{eqnarray*}
where
$
\Xi:=\frac{1}{n} \sum_{j=1}^n \xi_j E_{X_j,X_j'},  \xi_j:=
Y_j-S_{\ast}(X_j,X_j').
$
We can now rewrite the last bound as
\begin{eqnarray*}
&&
2\langle\hat S-S_{\ast}, \hat S-S\rangle_{L_2(P)} +
\eps\langle\hat{V},\hat{S}-S \rangle+ 2\bar\eps \bigl\langle W^{1/2}(
\hat{S}-S),W^{1/2}(\hat{S}-S) \bigr\rangle_{L_2(\Pi^2)}
\\
&&\qquad
\leq -2\bar\eps \bigl\langle W^{1/2}S,W^{1/2}(
\hat{S}-S) \bigr\rangle_{L_2(\Pi^2)} +2\langle\Xi,\hat{S}-S
\rangle\\
&&\quad\qquad{}+2(P-P_n) \bigl((\hat S-S_{\ast}) (\hat S-S)\bigr)
\end{eqnarray*}
and use a simple identity
\begin{eqnarray*}
2\langle\hat S-S_{\ast}, \hat S-S\rangle_{L_2(P)}&=&2\langle\hat
S-S_{\ast}, \hat S-S\rangle_{L_2(\Pi^2)}\\
&=& \|\hat{S}-S_\ast
\|^2_{L_2(\Pi^2)}+\|\hat S-S\|_{L_2(\Pi^2)}^2-
\|S-S_{\ast}\|_{L_2(\Pi^2)}^2
\end{eqnarray*}
to get the following bound:
%
\begin{eqnarray}
\label{eqstep1}&&
\|\hat{S}-S_\ast\|^2_{L_2(\Pi^2)}+
\|\hat S-S\|_{L_2(\Pi^2)}^2\nonumber\\
&&\quad{}+2\bar \eps\bigl\|W^{1/2}(\hat{S}-S)
\bigr\|_{L_2(\Pi^2)}^2
+\eps\langle\hat{V},\hat{S}-S \rangle
\nonumber\\
&&\qquad \leq\|S-S_{\ast}
\|_{L_2(\Pi^2)}^2-2\bar\eps\bigl\langle W^{1/2}S,W^{1/2}(
\hat{S}-S)\bigr\rangle_{L_2(\Pi^2)}
\\
&&\qquad\quad{}+
2\langle\Xi,\hat{S}-S\rangle + 2(P-P_n)
(S-S_{\ast}) (\hat S-S)\nonumber\\
&&\qquad\quad{} +2(P-P_n) (\hat S-S)^2.\nonumber
\end{eqnarray}

For an arbitrary $V \in\partial\| S\|_1$, $V=\operatorname{ sign}(S) +\mathcal
{P}_L^{\perp}(M)$, where $M$ is a matrix with $\|M\| \leq1$. It
follows from the trace duality property that there exists an $M$ with
$\| M \| \leq1$ [to be specific, $M=\operatorname{ sign}(\mathcal{P}_L^{\perp
}(\hat S))$]
such that
\[
\bigl\langle\mathcal{P}_L^{\perp}(M), \hat{S}-S \bigr
\rangle=\bigl\langle M,\mathcal{P}_L^{\perp}(\hat{S}-S) \bigr
\rangle=\bigl\langle M,\mathcal {P}_L^{\perp}(\hat S) \bigr
\rangle=\bigl\| \mathcal{P}_L^{\perp}(\hat{S})\bigr\|_1,
\]
where the first equality is based on the fact that $\mathcal
{P}_L^\perp$ is a self-adjoint operator and the second equality is
based on the fact that $S$ has support $L$.
Using this equation and monotonicity of\vadjust{\goodbreak} subdifferentials of convex
functions, we get
$
\langle\operatorname{ sign} (S),\hat{S}-S \rangle+ \| \mathcal{P}_L^{\perp
}(\hat{S})\|_1=\langle V, \hat{S}-S \rangle\leq\langle\hat{V},
\hat{S}-S \rangle.
$
Substituting this into the left-hand side of (\ref{eqstep1}), it is
easy to get
%
\begin{eqnarray}
\label{basic} && \|\hat S-S_{\ast}\|_{L_2(\Pi^2)}^2+ \|
\hat S-S\|_{L_2(\Pi^2)}^2\nonumber\\
&&\quad{}+\eps\bigl\|\mathcal{P}_L^{\perp}(\hat S)
\bigr\|_1 +2\bar\eps\bigl\|W^{1/2}(\hat S-S)\bigr\|_{L_2(\Pi^2)}^2
\nonumber\\
&&
\qquad
\leq\|S-S_{\ast}\|_{L_2(\Pi^2)}^2 -\eps\bigl
\langle\operatorname{ sign}(S),\hat S-S\bigr\rangle
\nonumber
\\[-8pt]
\\[-8pt]
\nonumber
&&\qquad\quad{}-2\bar\eps\bigl\langle W^{1/2}S, W^{1/2}(\hat
S-S)\bigr\rangle_{L_2(\Pi^2)}
\\
&&\qquad\quad{}
+ 2 \langle\Xi,\hat S-S\rangle + 2(P-P_n)
(S-S_{\ast}) (\hat S-S)\nonumber \\
&&\qquad\quad{}+2(P-P_n) (\hat S-S)^2.\nonumber
\end{eqnarray}

We need to bound the right-hand side of (\ref{basic}).
We start with deriving a bound on $\langle\operatorname{ sign}(S),\hat
S-S\rangle$, expressed in terms of function $\varphi$.
Note that, for all $\lambda>0$,
\begin{eqnarray*}
\bigl\langle\operatorname{ sign}(S),\hat S-S\bigr\rangle&=& \sum
_{k=1}^m \bigl\langle\operatorname{
sign}(S)\phi_k, (\hat S-S)\phi_k\bigr\rangle\\
&=&
\sum_{\lambda_k\leq\lambda} \bigl\langle\operatorname{
sign}(S)\phi_k, (\hat S-S)\phi_k\bigr\rangle\\
&&{} + \sum
_{\lambda_k>\lambda} \biggl\langle\frac{\operatorname{ sign}(S)\phi
_k}{\sqrt{\lambda_k}}, \sqrt{
\lambda_k}(\hat S-S)\phi_k \biggr\rangle,
\end{eqnarray*}
which easily implies
%
\begin{eqnarray} \label{odin} &&\bigl|\bigl\langle\operatorname{ sign}(S),\hat S-S\bigr\rangle\bigr|\nonumber\\
&&\qquad\leq
\biggl(\sum_{\lambda_k\leq\lambda}\bigl \|\operatorname{
sign}(S)\phi_k\bigr\| ^2 \biggr)^{1/2} \biggl(\sum
_{\lambda_k\leq\lambda}\bigl \|(\hat S-S)\phi_k
\bigr\|^2 \biggr)^{1/2}
\nonumber
\\[-8pt]
\\[-8pt]
\nonumber
&&\qquad\quad{}+
\biggl(\sum_{\lambda_k>\lambda} \frac{\|\operatorname{ sign}(S)\phi_k\|^2}{\lambda_k}
\biggr)^{1/2} \biggl(\sum_{\lambda_k>\lambda}
\lambda_k\bigl\|(\hat S-S)\phi_k\bigr\| ^2
\biggr)^{1/2} \\
&&\qquad\leq
\biggl(\sum_{\lambda_k\leq\lambda} \|P_L
\phi_k\|^2 \biggr)^{1/2}\| \hat S-S
\|_2 + \biggl(\sum_{\lambda_k>\lambda}
\frac{\|P_L\phi_k\|^2}{\lambda_k} \biggr)^{1/2} \bigl\|W^{1/2}(\hat S-S)
\bigr\|_2.\nonumber\hspace*{-35pt}
\end{eqnarray}

We will now use the following elementary lemma.

\begin{lemma}
\label{boundsums}
Let $c_{\gamma}:=\frac{c+\gamma}{\gamma}$.
For all $\lambda>0$,
\[
\sum_{\lambda_k>\lambda} \frac{\|P_L\phi_k\|^2}{\lambda_k}\leq
c_{\gamma}\frac{\varphi(\lambda)}{\lambda} \quad\mbox{and}\quad \sum
_{\lambda_k>\lambda} \frac{1}{\lambda_k}\leq c_{\gamma}
\frac{\bar F(\lambda)}{\lambda}.\vadjust{\goodbreak}
\]
\end{lemma}

\begin{pf}
Denote $H_k:=\sum_{j=1}^l \|P_L\phi_j\|^2, k=1,\ldots, m$.
Suppose that $\lambda\in [\lambda_l,\lambda_{l+1}]$ for some
$l=k_0-1, \ldots, m-1$.
We will use the properties of functions $\varphi\in\Psi_{S,W}$ and
$\bar F$. In particular, recall that the functions $\frac{\varphi
(\lambda)}{\bar F(\lambda)}$ and $\frac{\bar F(\lambda)}{\lambda}$
are nonincreasing.
Using these properties and the condition that $\lambda_{k+1}\leq
c\lambda_k, k\geq k_0$ we get
\begin{eqnarray*}
\label{dva}
\sum_{\lambda_k>\lambda}
\frac{\|P_L\phi_k\|^2}{\lambda_k}&=& \sum_{k=l+1}^{m-1}
H_k \biggl(\frac{1}{\lambda_k}-\frac{1}{\lambda
_{k+1}} \biggr) +
\frac{H_m}{\lambda_m}-\frac{H_l}{\lambda_{l+1}}\\
&\leq&
\sum_{k=l+1}^{m-1} \varphi(
\lambda_k) \biggl(\frac{1}{\lambda_k}-\frac{1}{\lambda
_{k+1}} \biggr) +
\frac{\varphi(\lambda_m)}{\lambda_m}\\
& \leq& c\sum_{k=l+1}^{m-1}
\frac{\varphi(\lambda_{k+1})}{\lambda_{k+1}^2}(\lambda _{k+1}-\lambda_{k}) +
\frac{\varphi(\lambda_m)}{\lambda_m}\\
& \leq&
c\int_{\lambda}^{\infty}\frac{\varphi(s)}{s^2}
\,ds +\frac{\varphi(\lambda)}{\lambda} \leq c\int_{\lambda}^{\infty}
\frac{\varphi(s)}{\bar F(s)} \frac{\bar F(s)}{s^2}\,ds + \frac{\varphi(\lambda)}{\lambda}\\
&\leq&
c\frac{\varphi(\lambda)}{\bar F(\lambda)}\int_{\lambda}^{\infty}
\frac{\bar F(s)}{s^2}\,ds + \frac{\varphi(\lambda)}{\lambda}\leq \frac{c}{\gamma}
\frac{\varphi(\lambda)}{\bar F(\lambda)} \frac{\bar F(\lambda)}{\lambda} + \frac{\varphi(\lambda)}{\lambda} \\
&=&\frac{c+\gamma}{\gamma}
\frac{\varphi(\lambda)}{\lambda},
\end{eqnarray*}
which proves the first bound. To prove the second bound, replace
in the inequalities above $\|P_L \phi_k\|^2$ by $1$ and $\varphi
(\lambda)$
by ${\bar F(\lambda)}$. In the case when $\lambda\geq\lambda_m$, both
bounds are trivial since their left-hand sides are equal to zero.
\end{pf}

It follows from from (\ref{odin}) and the first bound of Lemma \ref
{boundsums} that
%
\begin{eqnarray}
\label{tri} &&\bigl|\bigl\langle\operatorname{ sign}(S),\hat
S-S\bigr\rangle\bigr|\nonumber\\
&&\qquad\leq
\sqrt{\varphi(\lambda)}\|\hat S-S\|_2 + \sqrt{c_{\gamma}
\frac{\varphi(\lambda)}{\lambda}} \bigl\|W^{1/2}(\hat S-S)\bigr\|_2
\\
&&\qquad=
m\sqrt{\varphi(\lambda)}\|\hat S-S\|_{L_2(\Pi^2)} + m
\sqrt{c_{\gamma}\frac{\varphi(\lambda)}{\lambda}} \bigl\|W^{1/2}(\hat S-S)
\bigr\|_{L_2(\Pi^2)}.\nonumber
\end{eqnarray}
This implies the following bound:
%
\begin{eqnarray}
\label{tri'} &&\eps\bigl|\bigl\langle\operatorname{ sign}(S),\hat S-S\bigr\rangle\bigr|
\nonumber
\\
&&\qquad
\leq \varphi(\lambda)m^2\eps^2+ \frac{1}{4}\|\hat
S-S\|_{L_2(\Pi^2)}^2 +
c_{\gamma}\frac{\varphi(\lambda)}{\lambda}\frac{m^2\eps^2}{\bar
\eps} \\
&&\qquad\quad{}+
\frac{\bar\eps}{4}\bigl\|W^{1/2}(\hat
S-S)\bigr\|_{L_2(\Pi^2)}^2,\nonumber
\end{eqnarray}
where we used twice an elementary inequality $ab\leq a^2+\frac
{1}{4}b^2, a,b>0$.
We will apply this bound for $\lambda=\bar\eps^{-1}$ to get the
following inequality:
%
\begin{eqnarray}
\label{tri''} && \eps\bigl|\bigl\langle\operatorname{ sign}(S),\hat S-S\bigr\rangle\bigr|
\nonumber
\\
&&\qquad\leq
(c_{\gamma}+1)\varphi\bigl(\bar\eps^{-1}
\bigr)m^2\eps^2+ \frac{1}{4}\|\hat S-S
\|_{L_2(\Pi^2)}^2 \\
&&\qquad\quad{}+\frac{\bar\eps}{4}\bigl\|W^{1/2}(\hat S-S)
\bigr\|_{L_2(\Pi^2)}^2.\nonumber
\end{eqnarray}

To bound the next term in the right-hand side of (\ref{basic}),
note that
%
\begin{eqnarray}
\label{chetyre}&& \bar\eps\bigl|\bigl\langle W^{1/2}S, W^{1/2}(\hat
S-S)\bigr\rangle_{L_2(\Pi
^2)}\bigr|
\nonumber
\\[-8pt]
\\[-8pt]
\nonumber
&&\qquad\leq \bar\eps\bigl\|W^{1/2}S
\bigr\|_{L_2(\Pi^2)}^2 + \frac{\bar\eps}{4}\bigl\| W^{1/2}(\hat
S-S)\bigr\|_{L_2(\Pi^2)}^2.
\end{eqnarray}

The main part of the proof deals with bounding the stochastic term
\[
2\langle\Xi,\hat{S}-S\rangle + 2(P-P_n) (S-S_{\ast}) (\hat
S-S) +2(P-P_n) (\hat S-S)^2
\]
on the right-hand side of (\ref{basic}). To this end, define
(for fixed $S,S_{\ast}$)
\begin{eqnarray*}
f_{A}(y,u,v)&:=&\bigl(y-S_{\ast}(u,v)\bigr) (A-S)
(u,v)-(S-S_{\ast}) (u,v) (A-S) (u,v)
\\
&&{}- (A-S)^2(u,v)\\
&=& \bigl(y-S(u,v)\bigr) (A-S)
(u,v)-(A-S)^2(u,v),
\end{eqnarray*}
and consider the following empirical process:
\[
\alpha_n(\delta_1,\delta_2,
\delta_3):= \sup \bigl\{\bigl|(P_n-P) (f_{A})\bigr|
\dvtx A\in{\cal T}(\delta_1,\delta_2,
\delta_3) \bigr\},
\]
where
\begin{eqnarray*}
&\hspace*{-4pt}&{\cal T}(\delta_1,\delta_2,\delta_3)\\
&&\hspace*{-4pt}\qquad:=
\bigl\{A\in{\mathbb D}\dvtx \| A-S\|_{L_2(\Pi^2)}\leq\delta_1, \bigl\|{
\cal P}_{L}^{\perp}A\bigr\|_1\leq \delta_2,
\bigl\|W^{1/2}(A-S)\bigr\|_{L_2(\Pi^2)}\leq\delta_3 \bigr\}.
\end{eqnarray*}
Clearly, we have
%
\begin{eqnarray}
\label{stochterm} && 2\langle\Xi,\hat{S}-S \rangle + 2(P-P_n)
(S-S_{\ast}) (\hat S-S) +2(P-P_n) (\hat S-S)^2
\nonumber
\\[-8pt]
\\[-8pt]
\nonumber
&&\qquad\leq
2\alpha_n \bigl(\|\hat S-S\|_{L_2(\Pi^2)}, \bigl\|{\cal
P}_{L}^{\perp}\hat S\bigr\|_1,\bigl \|W^{1/2}(\hat
S-S)\bigr\|_{L_2(\Pi^2)} \bigr),
\end{eqnarray}
and it remains to provide an upper bound on $\alpha_n (\delta
_1,\delta_2, \delta_3)$
that is uniform in some intervals of the parameters $\delta_1, \delta
_2, \delta_3$
(such that either the norms $\|\hat S-S\|_{L_2(\Pi^2)},
\|{\cal P}_{L}^{\perp}\hat S\|_1, \|W^{1/2}(\hat S-S)\|_{L_2(\Pi^2)}$
belong to these intervals with a high probability, or bound of the theorem
trivially holds). Note that the functions $f_A$ are uniformly bounded by
a numerical constant (under the assumptions that $a=1$, $|Y|\leq a$ and
all the
kernels are also bounded by $a$) and we have
$
Pf_A^2 \leq c_1\|A-S\|_{L_2(\Pi)}^2
$
with some numerical constant $c_1>0$. Using Talagrand's concentration\vadjust{\goodbreak}
inequality for empirical processes we conclude that for fixed $\delta
_1, \delta_2, \delta_3$ with probability at least $1-e^{-t}$ and with
some constant \mbox{$c_2>0$}
$\alpha_n(\delta_1,\delta_2,\delta_3)\leq2 {\mathbb E}\alpha
_n(\delta_1,\delta_2, \delta_3)+ c_2 (\delta_1\sqrt{\frac
{t}{n}}+ \frac{t}{n} ).
$
We will make this bound uniform in $\delta_k\in[\delta_k^{-}, \delta
_k^+], \delta_k^-<\delta_k^+, k=1,2,3$ (these intervals will be
chosen later).
Define
$
\delta_k^j:=\delta_k^+ 2^{-j}, j=0,\ldots, [\log_2(\delta
_k^+/\delta_k^-)]+1,
k=1,2,3
$
and let
$
\bar t:= t+ \sum_{k=1}^3\log ([\log_2(\delta_k^+/\delta
_k^-)]+2 ).
$
By the union bound, with probability at least $1-e^{-t}$ and for all
$j_k=0, \ldots, [\log_2(\delta_k^+/\delta_k^-)]+1, k=1,2,3$,
$
\alpha_n(\delta_1^{j_1},\delta_2^{j_2},\delta_3^{j_3})\leq2
{\mathbb E}\alpha_n(\delta_1^{j_1},\delta_2^{j_2}, \delta_3^{j_3})+
c_2 (\delta_1^{j_1}\sqrt{\frac{\bar t}{n}}+ \frac{\bar
t}{n} ).
$
By monotonicity of $\alpha_n$ and of the right-hand side of the bound
with respect to each of the variables $\delta_1,\delta_2, \delta_3$,
we conclude that with the same probability and with some numerical
constant $c_3>0$,
for all $\delta_k \in[\delta_k^-,\delta_k^+], k=1,2,3$,
%
\begin{equation}
\label{talagr} \alpha_n(\delta_1,\delta_2,
\delta_3)\leq2 {\mathbb E}\alpha_n(2\delta_1,2
\delta_2, 2\delta_3)+ c_3 \biggl(
\delta_1\sqrt{\frac{\bar t}{n}}+ \frac{\bar t}{n} \biggr).
\end{equation}

To bound the expectation ${\mathbb E}\alpha_n(2\delta_1,2\delta_2,
2\delta_3)$
on the right-hand side of (\ref{talagr}), note that, by the definition
of function
$f_{A}$,
%
\begin{eqnarray}
\label{expal1} && {\mathbb E}\alpha_n(\delta_1,
\delta_2, \delta_3)\nonumber\\
&&\qquad\leq
{\mathbb E}\sup \bigl\{ \bigl|(P_n-P) (y-S) (A-S) \bigr|\dvtx A
\in{\cal T}(\delta_1,\delta_2, \delta_3)
\bigr\} \\[-2pt]
&&\qquad\quad{}+
{\mathbb E}\sup \bigl\{ \bigl|(P_n-P) (A-S)^2 \bigr|
\dvtx A\in{\cal T}(\delta_1,\delta_2,
\delta_3) \bigr\}.\nonumber
\end{eqnarray}
A standard application of symmetrization inequality followed by contraction
inequality for Rademacher sums (see, e.g.,~\cite{Ko27}, Chapter 2) yields
%
\begin{eqnarray}
\label{expal2} && {\mathbb E}\sup \bigl\{ \bigl|(P_n-P)
(A-S)^2 \bigr|\dvtx A\in{\cal T}(\delta_1,
\delta_2, \delta_3) \bigr\}
\nonumber
\\[-9pt]
\\[-9pt]
\nonumber
&&\qquad \leq
16 {\mathbb E}\sup \bigl\{ \bigl|R_n(A-S) \bigr|\dvtx A\in{\cal
T}(\delta_1,\delta_2, \delta_3) \bigr\}.
\end{eqnarray}
It easily follows from (\ref{expal1}) and (\ref{expal2}) that
%
\begin{eqnarray}
\label{expal3}  {\mathbb E}\alpha_n(\delta_1,
\delta_2, \delta_3)&\leq& {\mathbb E}\sup \bigl\{\bigl|\langle
\Xi_1, A-S\rangle \bigr|\dvtx A\in{\cal T}(\delta_1,
\delta_2, \delta_3) \bigr\}
\nonumber
\\[-9pt]
\\[-9pt]
\nonumber
&&{}+
16 {\mathbb E}\sup \bigl\{\bigl|\langle\Xi_2, A-S\rangle\bigr|
\dvtx A\in{\cal T}(\delta_1,\delta_2,
\delta_3) \bigr\},
\end{eqnarray}
where
$
\Xi_1:= \frac{1}{n}\sum_{j=1}^n
(Y_j-S(X_j,X_j'))E_{X_j,X_j'}-{\mathbb E}(Y-S(X,X'))E_{X,X'}
$
and
$
\Xi_2:= \frac{1}{n}\sum_{j=1}^n \eps_j E_{X_j,X_j'},
\{\eps_j\}$ being i.i.d. Rademacher random variables independent
of $(X_1,X_1',Y_1),\ldots, (X_n,X_n',Y_n)$.
We will upper bound the expectations on the right-hand side of (\ref
{expal3}),
which reduces to bounding
$
{\mathbb E}\sup \{|\langle\Xi_i, A-S\rangle|\dvtx
A\in{\cal T}(\delta_1,\delta_2, \delta_3) \}
$
for each of the random matrices $\Xi_1,\Xi_2$.
For $i=1,2$ and $A\in{\cal T}(\delta_1,\delta_2, \delta_3)$, we have
%
\begin{eqnarray}
\label{pjat}  \bigl| \langle\Xi_i, A-S \rangle\bigr|&\leq&\bigl |\bigl\langle
\Xi_i, {\cal P}_L (A-S)\bigr\rangle\bigr| + \bigl|\bigl\langle
\Xi_i, {\cal P}_L^{\perp} (A)\bigr\rangle\bigr|
\nonumber
\\[-2pt]
&\leq&\bigl|\langle{\cal P}_L\Xi_i, A-S\rangle\bigr|+
\|\Xi_i\|\bigl\|{\cal P}_L^{\perp}(A)
\bigr\|_1\\[-2pt] 
&\leq&\bigl|\langle{\cal P}_L
\Xi_i, A-S\rangle\bigr|+ \delta_2 \|\Xi_i\|.\nonumber
\end{eqnarray}
To bound $\|\Xi_i\|$, we use the following simple corollary of a
well-known noncommutative Bernstein inequality\vadjust{\goodbreak} (see, e.g., \cite
{Tropp}) obtained by integrating exponential tails of this inequality:
let $Z$ be a random symmetric matrix with $\EE Z=0$, $\sigma_Z^2:=\|
\EE Z^2\|$ and $\|Z\| \leq U$ for some $U > 0$ and let $Z_1,\ldots,Z_n$
be $n$ i.i.d. copies of $Z$. Then
%
\begin{equation}
\label{OpBer} {\mathbb E}\Biggl\llVert \frac{1}{n}\sum
_{j=1}^{n}Z_j \Biggr\rrVert \leq4
\biggl( \sigma_Z \sqrt{\frac{\log(2m)}{n}} \vee U
\frac{\log
(2m)}{n} \biggr).
\end{equation}
It is applied
to i.i.d. random matrices
\[
Z_j:=\bigl(Y_j -S\bigl(X_j,X_j'
\bigr)\bigr)E_{X_j,X_j'}-\EE\bigl(Y -S\bigl(X,X'
\bigr)\bigr)E_{X, X'}
\]
in the case of matrix $\Xi_1$ and to i.i.d. random matrices
$Z_j:=\eps_j E_{X_j,X_j'}$ in the case of matrix $\Xi_2$.
In both cases, $\|Z_j\|\leq4$ and, by a simple computation, $\sigma
_{Z_j}^2:=\|{\mathbb E}Z_j^2\|\leq4/m$ (see, e.g.,
\cite{Ko27}, Section 9.4),
bound (\ref{OpBer}) implies that, for $i=1,2$,
%
\begin{equation}
\label{expectXi} {\mathbb E}\| \Xi_i \|\leq16 \biggl[\sqrt{
\frac{\log(2m)}{n
m}}\vee \frac{\log(2m)}{n} \biggr]=: \eps^{\ast}.
\end{equation}

To control the term $|\langle{\cal P}_L\Xi_i, A-S\rangle|$ in bound
(\ref{pjat}),
we will use the following lemma.

\begin{lemma}
\label{lemma1}
For all $\delta>0$,
\[
{\mathbb E}\sup_{\|M\|_2\leq\delta, \|W^{1/2}M\|_2\leq1}\bigl|\langle {\cal
P}_L\Xi_i,M\rangle\bigr|\leq 4\sqrt{2}\sqrt{c_{\gamma}+1}
\sqrt{\frac{1}{nm}}\delta\sqrt {\varphi\bigl(\delta^{-2}\bigr)}.
\]
\end{lemma}

\begin{pf}
For all symmetric $m\times m$ matrices $M$,
\[
\langle{\cal P}_L \Xi_i, M\rangle= \sum
_{k,j=1}^m \langle{\cal P}_L
\Xi_i, \phi_k \otimes\phi _j\rangle \langle
M, \phi_k \otimes\phi_j\rangle.
\]
Assuming that
\begin{eqnarray*}
\|M\|_2^2&=&\sum_{k,j=1}^m
\bigl|\langle M, \phi_k\otimes\phi_j\rangle \bigr|^2
\leq \delta^2 \quad\mbox{and}\\[-2pt]
 \bigl\|W^{1/2}M\bigr\|_2^2
&=& \sum_{k,j=1}^m \lambda_k \bigl|
\langle M, \phi_k\otimes\phi_j\rangle \bigr|^2
\leq1,
\end{eqnarray*}
it is easy to conclude that
$
\sum_{k,j=1}^m \frac{|\langle M, \phi_k\otimes\phi_j\rangle
|^2}{\lambda_k^{-1}\wedge\delta^2}\leq2.
$
It follows
%
\begin{eqnarray}
\label{abab} && \qquad\bigl|\langle{\cal P}_L\Xi_i, M\rangle\bigr|
\nonumber\\[-2pt]
&&\qquad\leq
\Biggl(\sum_{k,j=1}^m \bigl(
\lambda_k^{-1}\wedge\delta^2\bigr)\bigl|\langle{\cal
P}_L \Xi, \phi_k \otimes\phi_j
\rangle\bigr|^2 \Biggr)^{1/2} \Biggl(\sum
_{k,j=1}^m \frac{|\langle M, \phi_k\otimes\phi_j\rangle|^2}{\lambda
_k^{-1}\wedge\delta^2} \Biggr)^{1/2}
\\
&&\qquad\leq
\sqrt{2} \Biggl(\sum_{k,j=1}^m
\bigl(\lambda_k^{-1}\wedge\delta^2\bigr)\bigl|
\langle{\cal P}_L \Xi, \phi_k \otimes\phi_j
\rangle\bigr|^2 \Biggr)^{1/2}.\nonumber
\end{eqnarray}
Consider the following inner product:
\[
\langle M_1,M_2 \rangle_{w}:= \sum
_{k,j=1}^m \bigl(\lambda_k^{-1}
\wedge \delta^2\bigr) \langle M_1, \phi_k
\otimes\phi_j\rangle \langle M_2, \phi_k
\otimes\phi_j\rangle,
\]
and let $\|\cdot\|_w$ be the corresponding norm. We will provide an upper
bound on
$
{\mathbb E}\|{\cal P}_L \Xi_i\|_w={\mathbb E}
(\sum_{k,j=1}^m (\lambda_k^{-1}\wedge\delta^2)|\langle{\cal
P}_L \Xi, \phi_k \otimes\phi_j\rangle|^2 )^{1/2}.
$
Recall that
$
\Xi_i = n^{-1}\sum_{j=1}^n \zeta_j E_{X_j,X_j'}- {\mathbb E}(\zeta E_{X,X'}),
$
where $\zeta_j=Y_j - S(X_j,X_j')$ for $i=1$ and $\zeta_j =\eps_j$
for $i=2$.
Note that in the first case $|\zeta_j|\leq2$, and in the second case
$|\zeta_j|\leq1$.
Therefore,
%
\begin{eqnarray}
\label{abub} {\mathbb E}\|{\cal P}_L \Xi_i
\|_w&\leq& {\mathbb E}^{1/2}\|{\cal P}_L
\Xi_i\|_w^2\leq \sqrt{
\frac{{\mathbb E}\zeta^2 \|{\cal P}_L E_{X,X'}\|_w^2}{n}}
\nonumber
\\[-9pt]
\\[-9pt]
\nonumber
& \leq&2\sqrt{\frac{{\mathbb E}\|{\cal P}_L E_{X,X'}\|_w^2}{n}}.
\end{eqnarray}
It remains to bound ${\mathbb E}\|{\cal P}_L E_{X,X'}\|_w^2$,
%
\begin{eqnarray}
\label{var} && {\mathbb E}\bigl\|{\cal P}_L (E_{X,X'})
\bigr\|_w^2\nonumber\\[-2pt]
&&\qquad=
{\mathbb E}\sum_{k,j=1}^m
\bigl(\lambda_k^{-1}\wedge\delta^2\bigr)\bigl|\bigl
\langle {\cal P}_L (E_{X,X'}), \phi_k \otimes
\phi_j\bigr\rangle\bigr|^2
\nonumber\\[-2pt]
&&\qquad=
\sum_{k,j=1}^m \bigl(
\lambda_k^{-1}\wedge\delta^2\bigr)
m^{-2}\sum_{u,v\in V}\bigl|\bigl\langle
E_{u,v}, {\cal P}_L(\phi_k \otimes\phi
_j)\bigr\rangle\bigr|^2
\nonumber\\[-2pt]
&&\qquad\leq
m^{-2}\sum_{k,j=1}^m
\bigl(\lambda_k^{-1}\wedge\delta^2\bigr) \bigl\|{
\cal P}_L(\phi_k \otimes\phi_j)
\bigr\|_2^2
\\[-2pt]
&&\qquad\leq
2m^{-2}\sum_{k,j=1}^m
\bigl(\lambda_k^{-1}\wedge\delta^2\bigr)
\bigl(\|P_L\phi_k\|^2+\|P_L
\phi_j\|^2\bigr)
\nonumber\\[-2pt]
&&\qquad=
2m^{-1}\sum_{k=1}^m
\bigl(\lambda_k^{-1}\wedge\delta^2\bigr)
\|P_L\phi_k\|^2 + 2m^{-2}\sum
_{k=1}^m \bigl(\lambda_k^{-1}
\wedge\delta^2\bigr)\|P_L\|_2^2
\nonumber\\[-2pt]
&&\qquad=
2m^{-1}\sum_{k=1}^m
\bigl(\lambda_k^{-1}\wedge\delta^2\bigr)
\|P_L\phi_k\|^2 + 2m^{-2}r\sum
_{k=1}^m \bigl(\lambda_k^{-1}
\wedge\delta^2\bigr).\nonumber
\end{eqnarray}
Note that
%
\begin{equation}
\label{hi-hi} \qquad\sum_{k=1}^m \bigl(
\lambda_k^{-1}\wedge\delta^2\bigr)
\|P_L\phi_k\|^2 \leq \delta^2
\sum_{\lambda_k\leq\delta^{-2}} \|P_L\phi_k
\|^2 + \sum_{\lambda_k>\delta^{-2}}\lambda_k^{-1}
\|P_L\phi_k\|^2.
\end{equation}

Using the first bound of Lemma~\ref{boundsums},
we get from (\ref{hi-hi}) that
%
\begin{eqnarray}
\label{au}  \sum_{k=1}^m \bigl(
\lambda_k^{-1}\wedge\delta^2\bigr)
\|P_L\phi_k\|^2& \leq& \delta^2
\varphi\bigl(\delta^{-2}\bigr)+c_{\gamma}\delta^2
\varphi\bigl(\delta^{-2}\bigr)
\nonumber
\\[-8pt]
\\[-8pt]
\nonumber
&=&(c_{\gamma}+1)
\delta^2 \varphi\bigl(\delta^{-2}\bigr).
\end{eqnarray}
We also have
$
\sum_{k=1}^m (\lambda_k^{-1}\wedge\delta^2)
\leq
\sum_{\lambda_k\leq\delta^{-2}} \delta^2
+
\sum_{\lambda_k>\delta^{-2}}\lambda_k^{-1},
$
which, by the second bound of Lemma~\ref{boundsums}, implies that
%
\begin{equation}
\label{ua}\quad \sum_{k=1}^m \bigl(
\lambda_k^{-1}\wedge\delta^2\bigr) \leq
\delta^2 \bar F\bigl(\delta^{-2}\bigr)+ c_{\gamma}
\delta^2\bar F\bigl(\delta^{-2}\bigr) \leq(c_{\gamma}+1)
\delta^2 \bar F\bigl(\delta^{-2}\bigr).
\end{equation}
Using bounds (\ref{var}), (\ref{au}) and (\ref{ua}) and
the fact that $\varphi(\lambda)\geq\frac{r}{m}\bar F(\lambda)$, we get
%
\begin{eqnarray}
\label{varA}  {\mathbb E}\bigl\|{\cal P}_L (E_{X,X'})
\bigr\|_w^2 &\leq& 2m^{-1}(c_{\gamma}+1)
\delta^2\varphi\bigl(\delta^{-2}\bigr)+
2m^{-2}r(c_{\gamma}+1)\delta^2 \bar F
\bigl(\delta^{-2}\bigr)
\nonumber\hspace*{-35pt}
\\[-8pt]
\\[-8pt]
\nonumber
& \leq& 4m^{-1}(c_{\gamma}+1)
\delta^2\varphi\bigl(\delta^{-2}\bigr).
\end{eqnarray}

The proof follows from (\ref{abab}), (\ref{abub}) and (\ref{varA}).
\end{pf}

Let $\delta:= \frac{\delta_1}{\delta_3}$.
Using Lemma~\ref{lemma1}, we get
\begin{eqnarray*}
&&\hspace*{-4pt}
{\mathbb E}\sup \bigl\{\bigl|\langle{\cal P}_L
\Xi_i, A-S\rangle\bigr|\dvtx A\in {\cal T}(\delta_1,
\delta_2,\delta_3) \bigr\}\\
&&\hspace*{-4pt}\qquad\leq
{\mathbb E}\sup \bigl\{\bigl|\langle{\cal P}_L
\Xi_i, A-S\rangle\bigr|\dvtx \| A-S\|_{L_2(\Pi^2)}\leq
\delta_1, \bigl\|W^{1/2}(A-S)\bigr\|_{L_2(\Pi^2)}\leq
\delta_3 \bigr\}\\
&&\hspace*{-4pt}\qquad =
{\mathbb E}\sup \bigl\{\bigl|\langle{\cal P}_L
\Xi_i, A-S\rangle\bigr|\dvtx \| A-S\|_2\leq
\delta_1 m, \bigl\|W^{1/2}(A-S)\bigr\|_2\leq
\delta_3 m \bigr\}\\
&&\hspace*{-4pt}\qquad \leq
\delta_3 m {\mathbb E}\sup \bigl\{\bigl|\langle{\cal
P}_L \Xi_i, A-S\rangle\bigr|\dvtx \|A-S\|_2\leq
\delta,\bigl \|W^{1/2}(A-S)\bigr\|_{L_2(\Pi^2)}\leq 1 \bigr\} \\
&&\hspace*{-4pt}\qquad\leq
4\sqrt{2}\delta_3 m \sqrt{c_{\gamma}+1} \sqrt{
\frac{1}{nm}}\delta\sqrt{\varphi\bigl(\delta^{-2}\bigr)}= 4
\sqrt{2}\sqrt{c_{\gamma}+1} \sqrt{\frac{m}{n}}\delta_1
\sqrt {\varphi\bigl(\delta^{-2}\bigr)}.
\end{eqnarray*}
In the case when $\delta^2\geq\bar\eps$, we get
\[
{\mathbb E}\sup \bigl\{\bigl|\langle{\cal P}_L
\Xi_i, A-S\rangle\bigr|\dvtx A\in {\cal T}(\delta_1,
\delta_2,\delta_3) \bigr\}\leq 4\sqrt{2}
\sqrt{c_{\gamma}+1} \delta_1\sqrt{\frac{m\varphi(\bar
\eps^{-1})}{n}}.
\]
In the opposite case, when $\delta^2< \bar\eps$, we use the fact that
the function $\frac{\varphi(\lambda)}{\lambda}=\frac{\varphi
(\lambda)}{\bar F(\lambda)} \frac{\bar F(\lambda)}{\lambda}$ is
nonincreasing.
This implies that
$\delta^2\varphi(\delta^{-2})\leq\bar\eps\varphi(\bar\eps^{-1})$,
and we get
\begin{eqnarray*}
&&
{\mathbb E}\sup \bigl\{\bigl|\langle{\cal P}_L
\Xi_i, A-S\rangle\bigr|\dvtx A\in {\cal T}(\delta_1,
\delta_2,\delta_3) \bigr\}\\
&&\qquad\leq
4\sqrt{2} \sqrt{c_{\gamma}+1} \sqrt{\frac{m}{n}}
\delta_1 \sqrt{\varphi\bigl(\delta^{-2}\bigr)} = 4\sqrt{2}
\sqrt{c_{\gamma}+1}\sqrt{\frac{m}{n}} \delta_3 \sqrt{
\delta^2\varphi\bigl(\delta^{-2}\bigr)}\\
&&\qquad \leq
4\sqrt{2} \sqrt{c_{\gamma}+1}\sqrt{\frac{m}{n}}
\delta_3 \sqrt{\bar\eps\varphi\bigl(\bar\eps^{-1}\bigr)} =
4\sqrt{2}\sqrt{c_{\gamma}+1} {\sqrt{\bar\eps}}\delta_3 \sqrt{
\frac{m\varphi(\bar\eps^{-1})}{n}}.
\end{eqnarray*}
We can conclude that
\begin{eqnarray*}
&&
{\mathbb E}\sup \bigl\{\bigl|\langle{\cal P}_L
\Xi_i, A-S\rangle\bigr|\dvtx A\in {\cal T}(\delta_1,
\delta_2,\delta_3) \bigr\}
\\
&&\qquad\leq
4\sqrt{2} \sqrt{c_{\gamma}+1} \delta_1\sqrt{
\frac{m\varphi(\bar
\eps^{-1})}{n}} + 4\sqrt{2}\sqrt{c_{\gamma}+1} {\sqrt{\bar\eps}}
\delta_3 \sqrt{\frac{m\varphi(\bar\eps^{-1})}{n}}.
\end{eqnarray*}
This bound will be combined with (\ref{pjat}) and (\ref{expectXi}) to
get that, for $i=1,2$,
\begin{eqnarray*}
&&
{\mathbb E}\sup \bigl\{\bigl|\langle\Xi_i, A-S\rangle\bigr|\dvtx
A\in{\cal T}(\delta_1,\delta_2,\delta_3)
\bigr\}\\
&&\qquad\leq \eps^{\ast}\delta_2 +
4\sqrt{2}\sqrt{c_{\gamma}+1} \delta_1\sqrt{
\frac{m\varphi(\bar
\eps^{-1})}{n}} + 4\sqrt{2}\sqrt{c_{\gamma}+1} {\sqrt{\bar\eps}}
\delta_3 \sqrt{\frac{m\varphi(\bar\eps^{-1})}{n}}.
\end{eqnarray*}
In view of (\ref{expal3}), this yields the bound
\[
{\mathbb E}\alpha_n(\delta_1,
\delta_2,\delta_3)\leq C'
\eps^{\ast}\delta_2 + C'\delta_1
\sqrt{\frac{m\varphi(\bar\eps^{-1})}{n}} + C'{\sqrt{\bar\eps}}
\delta_3 \sqrt{\frac{m\varphi(\bar\eps^{-1})}{n}}
\]
that holds with some constant $C'>0$ for all $\delta_1, \delta_2,
\delta_3>0$.
Using (\ref{talagr}), we conclude that for some constants $C$ and
for all $\delta_k \in[\delta_k^-,\delta_k^+], k=1,2,3$,
\[
\alpha_n(\delta_1,\delta_2,
\delta_3)\leq C \biggl[\delta_1\sqrt{
\frac{m\varphi(\bar\eps^{-1})}{n}} +\delta_1\sqrt{\frac{\bar t}{n}} +
\frac{\bar t}{n}+ \eps^{\ast}\delta_2 + {\sqrt{\bar\eps}}
\delta_3 \sqrt{\frac{m\varphi(\bar\eps^{-1})}{n}} \biggr]
\]
that holds with probability at least $1-e^{-t}$.
This yields the following upper bound on the stochastic term in (\ref{basic})
[see also (\ref{stochterm})]:
%
\begin{eqnarray}
\label{stochtermbound}
&& 2\langle\Xi,\hat{S}-S\rangle + 2(P-P_n)
(S-S_{\ast}) (\hat S-S) +2(P-P_n) (\hat S-S)^2
\nonumber\\
&&\qquad\leq
2C \biggl[\|\hat S-S\|_{L_2(\Pi^2)}\sqrt{\frac{m\varphi(\bar\eps
^{-1})}{n}} +
\|\hat S-S\|_{L_2(\Pi^2)}\sqrt{\frac{\bar t}{n}} + \frac{\bar t}{n}
\\
&&\hspace*{15pt}\quad\qquad{}+
\eps^{\ast}\|{\cal P}_L \hat S
\|_1 + {\sqrt{\bar\eps}}\bigl\|W^{1/2}(\hat S-S)
\bigr\|_{L_2(\Pi^2)} \sqrt{\frac{m\varphi(\bar\eps^{-1})}{n}} \biggr]\nonumber
\end{eqnarray}
that holds provided that
%
\begin{eqnarray}
\label{conddel} \|\hat S-S\|_{L_2(\Pi^2)}&\in&\bigl[\delta_1^-,
\delta_1^+\bigr],\qquad \bigl\|{\cal P}_L^{\perp} \hat S
\bigr\|_1\in\bigl[\delta_2^-,\delta_2^+\bigr],
\nonumber
\\[-8pt]
\\[-8pt]
\nonumber
\bigl\|W^{1/2}(\hat S-S)\bigr\|_{L_2(\Pi^2)}&\in&\bigl[\delta_3^-,
\delta_3^+\bigr].
\end{eqnarray}
We substitute bound (\ref{stochtermbound}) in (\ref{basic})
and further bound some of its terms as follows:
\begin{eqnarray*}
2C\|\hat S-S\|_{L_2(\Pi^2)}\sqrt{\frac{m\varphi(\bar\eps
^{-1})}{n}}&\leq& \frac{1}{8}
\|\hat S-S\|_{L_2(\Pi^2)}^2+ 8C^2 \frac{m\varphi(\bar
\eps^{-1})}{n},
\\
2C\|\hat S-S\|_{L_2(\Pi^2)}\sqrt{\frac{\bar t}{n}}&\leq& \frac{1}{8}
\|\hat S-S\|_{L_2(\Pi^2)}^2+ 8C^2 \frac{\bar t}{n}
\end{eqnarray*}
and
\begin{eqnarray*}
&& 2C\sqrt{\bar\eps}\bigl\|W^{1/2}(\hat S-S)\bigr\|_{L_2(\Pi^2)}\sqrt{
\frac
{m\varphi(\bar\eps^{-1})}{n}}\\
&&\qquad\leq \frac{1}{4}\bar\eps\bigl\|W^{1/2}(\hat S-S)
\bigr\|_{L_2(\Pi^2)}^2+ 4C^2 \frac
{m\varphi(\bar\eps^{-1})}{n}.
\end{eqnarray*}
We will also use (\ref{tri''}) to control the term $\eps|\langle\operatorname{
sign}(S), \hat S-S\rangle|$ in (\ref{basic}) and (\ref{chetyre}) to
control the term
$\bar\eps|\langle W^{1/2}S, W^{1/2}(\hat S-S)\rangle|$.
If condition (\ref{condeps}) holds with $D\geq32 C$, then
$\eps\geq2C \eps^{\ast}$.
By a simple algebra, it follows from (\ref{basic}) that
\begin{eqnarray*}
\|\hat S-S_{\ast}\|_{L_2(\Pi^2)}^2&\leq&
\|S-S_{\ast}\|_{L_2(\Pi^2)}^2+ C_1
m^2 \eps^2\varphi\bigl(\bar\eps^{-1}\bigr) +
C_1\frac{m\varphi(\bar\eps^{-1})}{n}\\
&&{} +\bar\eps\bigl\|W^{1/2}S
\bigr\|_{L_2(\Pi^2)}^2 + \frac{\bar t}{n}
\end{eqnarray*}
with some constant $C_1>0$. Since, under condition (\ref{condeps})
with $a=1$,
$m^2 \eps^2 \geq D^2 \frac{m\log(2m)}{n}\geq D^2 \frac{m}{n}$,
we can conclude that
%
\begin{eqnarray}
\label{finis} &&\|\hat S-S_{\ast}\|_{L_2(\Pi^2)}^2
\nonumber
\\[-8pt]
\\[-8pt]
\nonumber
&&\qquad\leq
\|S-S_{\ast}\|_{L_2(\Pi^2)}^2+ C_2
m^2 \eps^2\varphi\bigl(\bar\eps^{-1}\bigr) +
\bar\eps\bigl\|W^{1/2}S\bigr\|_{L_2(\Pi^2)}^2 + \frac{\bar t}{n}
\end{eqnarray}
with some constant $C_2>0$.

We still have to choose the values of $\delta_k^{-}, \delta_k^{+}$
and to handle the case when conditions (\ref{conddel}) do not hold.
First note that due to
the assumption that $\|S\|_{L_{\infty}}\leq1, S\in{\mathbb D}$,
we have
$\|\hat S-S\|_{L_2(\Pi)}\leq2$,
$\|{\cal P}_L^{\perp}\hat S\|_1\leq\|\hat S\|_1 \leq\sqrt{m}\|\hat
S\|_2
\leq m^{3/2}
$
and
$
\|W^{1/2}(\hat S-S)\|_{L_2(\Pi^2)}\leq2\sqrt{\lambda_m}.
$
Thus, we can set $\delta_1^+:=2, \delta_2^+:=m^{3/2}, \delta
_3^{+}:=2\sqrt{\lambda_m}$,
which guarantees that the upper bounds of (\ref{conddel}) are satisfied.
We will also set
$
\delta_1^{-}=\delta_2^{-}:= n^{-1/2}, \delta_3^{-}:=\sqrt{\frac
{\tilde\lambda}{n}}.
$
In the case when one of the lower bounds of (\ref{conddel}) does not hold,
we can still use inequality (\ref{stochtermbound}), but we have to replace
each of the norms $\|\hat S-S\|_{L_2(\Pi)}, \|{\cal P}_L^{\perp}\hat
S\|_1,
\|W^{1/2}(\hat S-S)\|_{L_2(\Pi^2)}$ which are smaller than the corresponding
$\delta_k^{-}$ by the quantity $\delta_k^{-}$. Then it is straightforward
to check that inequality (\ref{finis}) still holds for some value of constant
$C_2>0$.
With the above choice of $\delta_k^{-}, \delta_k^+$, we have
$
\bar t \leq
t+
3\log (2\log_2 n+\frac{1}{2}\log_2\frac{\lambda_m}{\tilde
\lambda}+2 )=t_{n,m}.
$
This completes the proof.
\end{pf*}




\printaddresses


\begin{thebibliography}{19}

\bibitem{Aubin}
\begin{bbook}[mr]
\bauthor{\bsnm{Aubin},~\bfnm{Jean-Pierre}\binits{J.-P.}} \AND
\bauthor{\bsnm{Ekeland},~\bfnm{Ivar}\binits{I.}}
(\byear{1984}).
\btitle{Applied Nonlinear Analysis}.
\bpublisher{Wiley}, \blocation{New York}.
\bid{mr={0749753}}
\bptok{imsref}%
\end{bbook}
\endbibitem

\bibitem{CandesPlan}
\begin{barticle}[mr]
\bauthor{\bsnm{Cand{\`e}s},~\bfnm{Emmanuel~J.}\binits{E.~J.}} \AND
\bauthor{\bsnm{Plan},~\bfnm{Yaniv}\binits{Y.}}
(\byear{2011}).
\btitle{Tight oracle inequalities for low-rank matrix recovery from a minimal
number of noisy random measurements}.
\bjournal{IEEE Trans. Inform. Theory}
\bvolume{57}
\bpages{2342--2359}.
\bid{doi={10.1109/TIT.2011.2111771}, issn={0018-9448}, mr={2809094}}
\bptok{imsref}%
\end{barticle}
\endbibitem

\bibitem{CandesRecht}
\begin{barticle}[mr]
\bauthor{\bsnm{Cand{\`e}s},~\bfnm{Emmanuel~J.}\binits{E.~J.}} \AND
\bauthor{\bsnm{Recht},~\bfnm{Benjamin}\binits{B.}}
(\byear{2009}).
\btitle{Exact matrix completion via convex optimization}.
\bjournal{Found. Comput. Math.}
\bvolume{9}
\bpages{717--772}.
\bid{doi={10.1007/s10208-009-9045-5}, issn={1615-3375}, mr={2565240}}
\bptok{imsref}%
\end{barticle}
\endbibitem

\bibitem{CandesTao}
\begin{barticle}[mr]
\bauthor{\bsnm{Cand{\`e}s},~\bfnm{Emmanuel~J.}\binits{E.~J.}} \AND
\bauthor{\bsnm{Tao},~\bfnm{Terence}\binits{T.}}
(\byear{2010}).
\btitle{The power of convex relaxation: Near-optimal matrix completion}.
\bjournal{IEEE Trans. Inform. Theory}
\bvolume{56}
\bpages{2053--2080}.
\bid{doi={10.1109/TIT.2010.2044061}, issn={0018-9448}, mr={2723472}}
\bptok{imsref}%
\end{barticle}
\endbibitem

\bibitem{Gine}
\begin{bbook}[mr]
\bauthor{\bparticle{de~la}
\bsnm{Pe{\~n}a},~\bfnm{V{\'{\i}}ctor~H.}\binits{V.~H.}} \AND
\bauthor{\bsnm{Gin{\'e}},~\bfnm{Evarist}\binits{E.}}
(\byear{1999}).
\btitle{Decoupling: From Dependence to Independence}.
\bpublisher{Springer}, \blocation{New York}.
\bid{doi={10.1007/978-1-4612-0537-1}, mr={1666908}}
\bptnote{check year}%
\bptok{imsref}%
\end{bbook}
\endbibitem

\bibitem{GaifasLecue}
\begin{barticle}[mr]
\bauthor{\bsnm{Ga{\"{\i}}ffas},~\bfnm{St{\'e}phane}\binits{S.}} \AND
\bauthor{\bsnm{Lecu{\'e}},~\bfnm{Guillaume}\binits{G.}}
(\byear{2011}).
\btitle{Hyper-sparse optimal aggregation}.
\bjournal{J. Mach. Learn. Res.}
\bvolume{12}
\bpages{1813--1833}.
\bid{issn={1532-4435}, mr={2819018}}
\bptok{imsref}%
\end{barticle}
\endbibitem

\bibitem{Gross-2}
\begin{barticle}[mr]
\bauthor{\bsnm{Gross},~\bfnm{David}\binits{D.}}
(\byear{2011}).
\btitle{Recovering low-rank matrices from few coefficients in any basis}.
\bjournal{IEEE Trans. Inform. Theory}
\bvolume{57}
\bpages{1548--1566}.
\bid{doi={10.1109/TIT.2011.2104999}, issn={0018-9448}, mr={2815834}}
\bptok{imsref}%
\end{barticle}
\endbibitem

\bibitem{Klopp}
\begin{bmisc}[auto:STB|2013/03/04|13:35:07]
\bauthor{\bsnm{Klopp},~\bfnm{O.}\binits{O.}}
(\byear{2012}).
\bhowpublished{Noisy low-rank matrix completion with general sampling
distribution. Preprint}.
\bptok{imsref}%
\end{bmisc}
\endbibitem

\bibitem{Ko27}
\begin{bbook}[mr]
\bauthor{\bsnm{Koltchinskii},~\bfnm{Vladimir}\binits{V.}}
(\byear{2011}).
\btitle{Oracle Inequalities in Empirical Risk Minimization and Sparse Recovery
Problems}.
\bseries{Lecture Notes in Math.}
\bvolume{2033}.
\bpublisher{Springer}, \blocation{Heidelberg}.
\bid{doi={10.1007/978-3-642-22147-7}, mr={2829871}}
\bptok{imsref}%
\end{bbook}
\endbibitem

\bibitem{Ko26}
\begin{barticle}[mr]
\bauthor{\bsnm{Koltchinskii},~\bfnm{Vladimir}\binits{V.}},
\bauthor{\bsnm{Lounici},~\bfnm{Karim}\binits{K.}} \AND
\bauthor{\bsnm{Tsybakov},~\bfnm{Alexandre~B.}\binits{A.~B.}}
(\byear{2011}).
\btitle{Nuclear-norm penalization and optimal rates for noisy low-rank matrix
completion}.
\bjournal{Ann. Statist.}
\bvolume{39}
\bpages{2302--2329}.
\bid{doi={10.1214/11-AOS894}, issn={0090-5364}, mr={2906869}}
\bptok{imsref}%
\end{barticle}
\endbibitem

\bibitem{KoRan}
\begin{bmisc}[auto:STB|2013/03/04|13:35:07]
\bauthor{\bsnm{Koltchinskii},~\bfnm{V.}\binits{V.}} \AND
\bauthor{\bsnm{Rangel},~\bfnm{P.}\binits{P.}}
(\byear{2012}).
\bhowpublished{Low rank estimation of similarities on graphs. Preprint}.
\bptok{imsref}%
\end{bmisc}
\endbibitem


\bibitem{Massart}
\begin{bbook}[mr]
\bauthor{\bsnm{Massart},~\bfnm{Pascal}\binits{P.}}
(\byear{2007}).
\btitle{Concentration Inequalities and Model Selection}.
\bseries{Lecture Notes in Math.}
\bvolume{1896}.
\bpublisher{Springer}, \blocation{Berlin}.
\bid{mr={2319879}}
\bptok{imsref}%
\end{bbook}
\endbibitem

\bibitem{Negahban}
\begin{barticle}[auto:STB|2013/03/04|13:35:07]
\bauthor{\bsnm{Negahban},~\bfnm{S.}\binits{S.}} \AND
\bauthor{\bsnm{Wainwright},~\bfnm{M.~J.}\binits{M.~J.}}
(\byear{2012}).
\btitle{Restricted strong convexity and weighted matrix
completion: optimal bounds with noise}.
\bjournal{J. Mach. Learn. Res.}
\bvolume{13}
\bpages{1665--1697}.
\bptok{imsref}%
\end{barticle}
\endbibitem

\bibitem{Recht1}
\begin{barticle}[mr]
\bauthor{\bsnm{Recht},~\bfnm{Benjamin}\binits{B.}},
\bauthor{\bsnm{Fazel},~\bfnm{Maryam}\binits{M.}} \AND
\bauthor{\bsnm{Parrilo},~\bfnm{Pablo~A.}\binits{P.~A.}}
(\byear{2010}).
\btitle{Guaranteed minimum-rank solutions of linear matrix equations via
nuclear norm minimization}.
\bjournal{SIAM Rev.}
\bvolume{52}
\bpages{471--501}.
\bid{doi={10.1137/070697835}, issn={0036-1445}, mr={2680543}}
\bptok{imsref}%
\end{barticle}
\endbibitem

\bibitem{Rohde}
\begin{barticle}[mr]
\bauthor{\bsnm{Rohde},~\bfnm{Angelika}\binits{A.}} \AND
\bauthor{\bsnm{Tsybakov},~\bfnm{Alexandre~B.}\binits{A.~B.}}
(\byear{2011}).
\btitle{Estimation of high-dimensional low-rank matrices}.
\bjournal{Ann. Statist.}
\bvolume{39}
\bpages{887--930}.
\bid{doi={10.1214/10-AOS860}, issn={0090-5364}, mr={2816342}}
\bptok{imsref}%
\end{barticle}
\endbibitem

\bibitem{Tropp}
\begin{barticle}[auto:STB|2013/03/04|13:35:07]
\bauthor{\bsnm{Tropp},~\bfnm{J.~A.}\binits{J.~A.}}
(\byear{2012}).
\btitle{User-friendly tail bounds for sums of random matrices}.
\bjournal{Found. Comput. Math.}
\bvolume{12}
\bpages{389--434}.
\bptok{imsref}%
\end{barticle}
\endbibitem

\bibitem{Tsybakov}
\begin{bbook}[mr]
\bauthor{\bsnm{Tsybakov},~\bfnm{Alexandre~B.}\binits{A.~B.}}
(\byear{2009}).
\btitle{Introduction to Nonparametric Estimation}.
\bpublisher{Springer}, \blocation{New York}.
\bid{doi={10.1007/b13794}, mr={2724359}}
\bptok{imsref}%
\end{bbook}
\endbibitem

\end{thebibliography}
\end{document}